\documentclass[12pt]{article}

\usepackage[psamsfonts]{amssymb}
\usepackage{amsfonts,euscript,wrapfig}
\usepackage{amsmath,amsthm,graphicx}
\usepackage{booktabs}

\title{Self-Similar Jordan Arcs Which Do Not Satisfy OSC.}

\author{Kirill Kamalutdinov \and 
Andrey Tetenov\footnote{Supported by Russian Foundation of Basic Research project
13-01-00513} \and Dmitry Vaulin}

\begin{document}

\maketitle

\newcommand{\rr}{\mathbb{R}}
\newcommand \nn {\mathbb{N}}
\newcommand \zz {\mathbb{Z}}
\newcommand \bbc {\mathbb{C}}
\newcommand \rd {\mathbb{R}^d}

 \newcommand {\al} {\alpha}
\newcommand {\be} {\beta}
\newcommand {\da} {\delta}
\newcommand {\Da} {\Delta}
\newcommand {\ga} {\gamma}
\newcommand {\ka} {\varkappa}
\newcommand {\Ga} {\Gamma}
\newcommand {\la} {\lambda}
\newcommand {\La} {\Lambda}
\newcommand{\om}{\omega}
\newcommand{\Om}{\Omega}
\newcommand {\sa} {\sigma}
\newcommand {\Sa} {\Sigma}
\newcommand {\te} {\theta}
\newcommand {\vte} {\vartheta}
\newcommand {\fy} {\varphi}
\newcommand {\Fy} { \fy}
\newcommand {\ep} {\varepsilon}
\newcommand{\e}{\varepsilon}
\newcommand{ \vro}{\varrho}

\newcommand{\VEC}{\overrightarrow}
\newcommand{\IN}{{\subset}}
\newcommand{\NI}{{\supset}}
\newcommand \dd  {\partial}
\newcommand {\mmm}{{\setminus}}
\newcommand{\probel}{\vspace{.5cm}}
\newcommand{\8}{{\infty}}
\newcommand{\0}{{\varnothing}}
\newcommand{\vse}{$\blacksquare$}

\newcommand {\bfep} {{{\bar \varepsilon}}}
\newcommand {\Dl} {\Delta}
\newcommand{\vA}{{\vec {A}}}
\newcommand{\vB}{{\vec {B}}}
\newcommand{\vF}{{\vec {F}}}
\newcommand{\vf}{{\vec {f}}}
\newcommand{\vh}{{\vec {h}}}
\newcommand{\vJ}{{\vec {J}}}
\newcommand{\vK}{{\vec {K}}}
\newcommand{\vP}{{\vec {P}}}
\newcommand{\vX}{{\vec {X}}}
\newcommand{\vY}{{\vec {Y}}}
\newcommand{\vZ}{{\vec {Z}}}
\newcommand{\vx}{{\vec {x}}}
\newcommand{\va}{{\vec {a}}}
\newcommand{\vga}{{\vec {\gamma}}}

\newcommand{\eS}{{\EuScript S}}
\newcommand{\eH}{{\EuScript H}}
\newcommand{\eC}{{\EuScript C}}
\newcommand{\eP}{{\EuScript P}}
\newcommand{\eT}{{\EuScript T}}
\newcommand{\eG}{{\EuScript G}}
\newcommand{\eK}{{\EuScript K}}
\newcommand{\eF}{{\EuScript F}}
\newcommand{\eZ}{{\EuScript Z}}
\newcommand{\eL}{{\EuScript L}}
\newcommand{\eD}{{\EuScript D}}
\newcommand{\E}{{\EuScript E}}
\def \diam {\mathop{\rm diam}\nolimits}
\def \fix {\mathop{\rm fix}\nolimits}
\def \Lip {\mathop{\rm Lip}\nolimits}

\newcommand{\io}{I^\infty}
\newcommand{\hio}{\hat{I}^\infty}
\newcommand{\imo}{I^{-\omega}}
\newcommand{\is}{I^*}

\newcommand{\Be}{{{\bf e}}}
\newcommand{\bi}{{{\bf i}}}
\newcommand{\bj}{{{\bf j}}}
\newcommand{\bk}{{{\bf k}}}
\newcommand{\bn}{{{\bf n}}}
\newcommand{\bx}{{{\bf x}}}
\newcommand{\by}{{{\bf y}}}

\newtheorem{thm}{\bf Theorem}
 \newtheorem{cor}[thm]{\bf Corollary}
 \newtheorem{lem}[thm]{\bf Lemma}
 \newtheorem{prop}[thm]{\bf Proposition}
 \newtheorem{dfn}[thm]{\bf Definition}

\newcommand{\dok}{{\bf{Proof}}}


The  problem of  finding  explicit  geometrical criteria for a system $\EuScript{S}$ of contraction similarities, implying the open set condition,  is  discussed since 80-ies and  still remains open.

One  of  such  criteria  is the   {\it
finite intersection property}:\\    The  system
$\EuScript{S}=\{S_1,...,S_{m}  \}$with the attractor
 $K$,  has f.i. property if  for  any  t $i\neq j$ the  set
$S_i(K)\cap S_j(K)$ is  finite.

\medskip

It was proved in
2007 by C.Bandt and  H.Rao \cite{BR} that if a system $\EuScript{S}=\{S_1,...,S_{m}  \}$  of contraction similarities in $\mathbb R^2$ with a connected attractor $K$ has the finite intersection property, then it satisfies OSC.  The  authors wrote they  believe that   in $\mathbb R^3$ this is not  so.

 In this paper  we prove the following
\begin{thm}\label{main}
There is such system $\eS=\{S_1,...,S_{m}  \}$ of contraction similarities in $\rr^3$, which:\\
(1) does not satisfy OSC, \\  (2) satisfies one-point intersection property and \\ (3) whose attractor is
a Jordan arc $\ga\IN\rr^3$.
\end{thm}

Notice  that  the  statements (1) and (2) imply  that the  system $\eS$ does  not  satisfy weak separation property WSP.

To show  the  existence of  such  self-similar   arcs  we use the zipper  construction \cite{ATK} and prove  the  following  three  statements, which are  the  foundation  of  our  approach to construction of non-WSP curves:

First is  a {\bf general position   theorem}  for   fractal
curves,  which gives  the  condition specifying how to get  rid of
intersection  by  small deformations of pairs of fractal curves
within a given family  of such  pairs
$\{(\varphi(x,t),\psi(x,t))\}$, depending on some parameter $x\in
\rr^3$.

\begin{thm}\label{genpos}
Let  $\varphi(x,t),\psi(x,t):B^3\times I\to \mathbb R^3$  be  continuous maps  which\\
 (1) are $\alpha$-H\"older  with  respect  to $t$
 and\\  (2) satisfy the  condition: for  any $t,s\in I$  the  function $f(x,t,s)=\varphi(x,t)-\psi(x,s)$ is bi-Lipschitz  with respect   to $x$. \\
Then Hausdorff dimension of the  set $\Delta=\{x\in B^3| \varphi(x,I)\cap\psi(x,I)\neq\0\} $ does  not  exceed  ${2}/{\alpha}$.
\end{thm}

Second  is a  {\bf corollary of Barnsley's Collage Theorem}, that gives  the conditions under which a  deformation of a  fractal  is  bi-Lipschitz on its certain subpieces.

\begin{prop}\label{collage}
Let $\eS=\{S_1,...,S_m  \}, \eT=\{T_1,...,T_m  \}$ be  systems of  contractions   in  a  complete metric  space $X$, and  $q=\max(\Lip S_i,\Lip T_i)$. \\
Let $\fy:\io\to K$, $\psi:\io\to L$ be   index  maps for  the attractors $K(\eS)$ and $L(\eT)$.
Let $V$ be such bounded set, that for  any $i=1,...,m$, $ S_i(V)\IN V$  and  $ T_i(V)\IN V$.
Put $\Da_i(x)= T_i(x)- S_i(x)$ and suppose $\|\Da_i(x)\|\le \da$ for  any $i$  and  any $x\in V$. Then,

{\bf B1)}      If for some multiindex $\bj$,  $\da_1<\|\Da_\bj(x)\|(<\da_2$ for  any $x\in V$ and $\da_1>\dfrac{q_\bj\da}{1-q}$,   then for $\sa\in \tau_\bj\io$,
$$\da_1-\dfrac{q_\bj\da}{1-q}\le \|\psi(\sa)-\fy(\sa)\|\le\da_2+\dfrac{q_\bj\da}{1-q}$$

{\bf B2)}      If for some multiindices $\bi,\   \bj$,  $\da_1<|\Da_i(x)-\Da_\bj(y)|<\da_2$ for  any $x,y\in V$ and $\da_1>\dfrac{(q_\bi+q_\bj)\da}{1-q}$, then for $\tau\in \tau_\bj\io$, $\sa\in \sa_\bi\io$,
$$\da_1-\dfrac{(q_\bi+q_\bj)\da}{1-q}\le  \|\psi(\sa)-\fy(\sa)-\psi(\tau)+\fy(\tau)\|\le\da_2+\dfrac{(q_\bi+q_\bj)\da}{1-q}$$
\end{prop}

And the third, most delicate, statement  allows to make  deformations  of system $\eS$ and  the  curve $\ga$, under which the violation of  WSP by the  system $\eS$ is  preserved.

\begin{thm}\label{2gen}
1) If  two-generator  subgroup $G=\langle \xi,\eta,\cdot\rangle$, $\xi=r e^{i\alpha},\eta=R e^{i\beta}$  in  $\bbc\mmm\{0\}$  is a dense subgroup  of  second type, then for  any   $z_1,z_2\in \mathbb{C}\setminus\{0\}$ there is such sequence $\{(n_k,m_k)\}$  that $\lim\limits_{k\to \infty}  \dfrac{z_1\xi^{n_k}}{z_2\eta^{m_k}}=1$,  $\lim\limits_{k\to\infty} e^{i n_k\alpha}=e^{-i \arg(z_1)}$,  $\lim\limits_{k\to\infty} e^{i n_k\beta}=e^{-i \arg(z_2)}$.\\
2) The set $\{(\xi,\eta)\}$ of all pairs of generators of dense subgroups of the second type, is dense in $\bbc^2$.
\end{thm}

The  subgroups mentioned  in the  theorem were  defined in \cite{TKV} and  we  give   a  short  summary  of  the  results  of  this  work in the   Addendum.

\probel

\section{ Brief  outline  of  the  method.}

The idea of the example is quite simple and  is  based  on zipper  construction introduced  by Vladislav Aseev  in \cite{ATK}.\\
\smallskip

Namely, a system $\eS=\{S_1,...,S_m\}$ of contractions  of  a metric  space $X$ is called a {\em zipper} with vertices $\{z_{0},...,z_{m}\} $ and signature $\e=(\e_{1},..,\e_{m}), \e_i\in\{0,1\}$, if for  any $i=1,...,m$,   $S_i(z_{0})=z_{i-1+\e_{i}}$ and $S_{i}(z_{m})=z_{i-\e_{i}}$.

We  denote the  attractor of  a zipper $\eS$ by $\ga(\eS)$, or  simply by $\ga$.

Note, that for  any zipper $\eS$ and any linear zipper $\eT$ on $[0,1]$  having   the  same  signature  $\e$  there  is unique  a  H\"older continuous $(\eS,\eT)$-equivariant map $\fy_{_{\eS\eT}}:[0,1]\to\ga$ 
 which  is  called  a  {\it linear  parametrizatiion}  of  $\eS$. \cite{ATK}

\medskip

 We begin with a self-similar zipper $\EuScript S=\{S_1,...,S_{2m}\}$ in $\mathbb R^3$ , whose  vertices $z_0,...,z_{2m}$ and  similarities $S_i$ are chosen in such way,  that  for some  closed bicone $V$ with vertices $z_0,z_{2m}$,

\medskip

{\bf(A1)} \ \   for any $S_i\in \EuScript S$,  $S_i(V)\subset V$  ;

\medskip

{\bf(A2)} \ \  for any such $i, j$, that $|j-i|>1$, \ \
 $S_i(V)\cap S_j(V)=\0$

\medskip

{\bf(A3)}\ \  for any $i\neq m+1$ and  $1\le i \le 2m$, \ \ \  $S_{i-1}(V)\cap S_i(V)=\{z_i\}$;

\medskip

At the same time, inside $S_m(V)\cup S_{m+1}(V)$ we provide that

{\bf(A4)} \ \  There   are such subarcs   $\ga_A\NI \ga_{m-3}$, $\ga_B\NI \ga_{m+4}$
and  such sequences $\{i_k\}, \{j_k\}$, that the  set  $S_m(\ga)\cap S_{m+1}(\ga)\mmm\{z_m\}$ is  a disjoint  union of sets

$$S_{m+1}S_{1}^{i_k}(\ga_A)\cap S_{m}S_{2m}^{j_k}(\ga_B)$$

\medskip

{\bf(A5)} \ \  The sequence of pairs of maps
 $S'_k=S_{m+1}S_{1}^{i_k},S''_k=S_{m}S_{2m}^{j_k}$, contains such subsequence $\{S'_{k_n}, S''_{k_n}\}$, that
$$\lim\limits_{n\to\infty}(S'_{k_n}S_{m-3})^{-1}(S''_{k_n}S_{m+4})={\rm Id}.$$

{\bf(A6)} \ \  There  is  such  linear zipper $\eT$ in $[0,1]$ that  the   H\"older  exponent  of  the linear parametrization $\fy_{_{\eS\eT}}$ is  greater  than  3/4.

\vspace{1cm}

The  condition (A5) means that  the system $\EuScript S$ does not satisfy WSP, but (A1) -- (A3) do not imply  in general, that $\gamma$ is a Jordan arc.

\bigskip

So, the main difficulty is to change the system $\EuScript S$ slightly in such way  that $\gamma$ gets rid of all its self-intersections without violating (A4, A5).

For that reason, instead of a single  zipper $\eS$, in Section 2 we  construct a  family  of self-similar zippers $\EuScript S_\xi=\{S_1,...,S_{2m}   \}$,  depending continuously on a parameter $\xi\in D$, where $D $ is  some domain in $\mathbb R^3$, so that for any $\xi\in D$, $\EuScript S_\xi$ satisfies (A1) - (A6)  with the same domain $V$, sequences $S'_k, S''_k$  and the same linear zipper $\eT$ for all $\xi\in D$.

The most delicate problem here is to choose such parameters for the  family $\{\eS_\xi\}$, that  the condition (A4) would  hold for all $\eS_\xi$. This is guaranteed by our Theorem 4 on the properties of dense subgroups in $\bbc^* $, proved in  \cite{TKV}.

After that we show  that  the  set of the parameters $\xi\in D$, for which $\EuScript S_\xi$ defines a Jordan arc, is dense in $D$.  To obtain this, we apply  the general position Theorem 2 the  following  way:

We take the family $\{\fy^\xi,\xi\in D\}$ of linear
parametrisations $\fy^\xi(t):[0,1]\to\ga^\xi$  of zippers
$\eS_\xi$ by the  zipper $\eT=\{T_1,...,T_m\}$. Denoting by  $\fy_k(\xi,t),
\psi_k(\xi,t)$  the parametrisations  of  the subarcs
$S'_k(\ga_A), S''_k(\ga_B)$, obtained by  restriction of
$\fy^\xi(t)$ to  the subintervals
$T_{m+1}T_1^{i_k}(I_A),T_{m}T_{2m}^{j_k} (I_B)$ of  $I=[0,1]$, we
consider the functions
$f_k(\xi,t,s)=\varphi_k(\xi,t)-\psi_k(\xi,s)$. Applying
Proposition 3, we  show that
 the  function $f_k(\xi,t,s)=\varphi_k(\xi,t)-\psi_k(\xi,s)$, is  bi-Lipschitz   with respect  to $\xi$   for  each fixed $t, s$. By our  construction,  $f_k(\xi,t,s)$ is  $\al$-H\"older  with respect  to t and s and $\al>2/3$.

Therefore, applying   Theorem \ref{genpos} to each  pair of  subarcs  from
the  sequence $S'_k(\ga_A),S''_k(\ga_B)$ and  using  Baire
category argument, we get  that each ball in the  domain  $D$
contains such $\xi$, that all the intersections
$S'_k(\ga_A)\cap S''_k(\ga_B)$ are  empty and therefore
 $\gamma^\xi$ is  a Jordan arc.
So the  set  of  parameters $\xi$, for  which the  system $\eS(\xi)$ has  a  Jordan attractor $\ga$ is  dense in $D$.

\section{ Setting  the  parameters of   zippers $\eS^\xi$.}

In this  section  we  define the parameter domain $D$, the  zippers $\eS^\xi$  for  any $\xi\in D$ and  the  bicones $V_i$ .

\subsection{ Polygonal lines  defining the  zippers $\eS^\xi$.}

First, we  define  three  angles $\be_1<\be_0<\be_2$, necessary  for  our  construction.\\
Take $\be_0=\arctan (1/2)$; let $\mu=0.0100512$ and put $\be_1=\be_0-2\mu$,
$\be_2=\be_0+\mu$.

\bigskip

We  define the domain $D$    in $\rr^3$  by  the  equation $$D=(1/1.02,1.02)\times (\be_0-\mu,\be_0-\mu/2)\times (-\mu,\mu)$$

Denote a point in $D$ by  $\xi=(\rho,\te,\fy)$.

Now  we  define a  family of zippers $\eS_\xi$ depending  on a  parameter $\xi\in D$

So  for  each $\xi\in D$ we  define a  polygonal line  in the  plane XY with   vertices $z_0
,...., z_{2m}$.   All its vertices, except $z_{m+1}$, do  not   depend  on  $\xi$.

The similarities $S_i$ in each  zipper  $\eS^\xi$ are the compositions  of  XY-plane preserving  similarities sending $\{z_0,z_{2m}\}$ to $\{z_{i-1},z_i\}$ and  rotations in  some  angle $\al_i$  around  the  axis $\{z_{i-1},z_i\}$.

In this  family, only the  point $z_{m+1}$ and   the  maps
$S_{m+1}, S_{m+2}$ and $S_{m+4}$ depend  on $\xi$, others  being
the  same  for  any $\xi\in D$.

Mentioning the  points, maps  or subset  of  the  attractor  of  each  system  $\eS_\xi$  we will write $z_i^\xi$ or $S_i^\xi$ only  in the  cases when it  is  needed for our  argument,
otherwise we  will  not  mention  the  parameter $\xi$, assuming the  dependence by  default.

\medskip

The  table below shows  the x and y   coordinates  of  the
vertices $z_i$,  the  contraction  ratios
$q_i=\dfrac{||z_i-z_{i-1}||}{||z_{2m}-z_0||}$ and rotation  angles $\al_i$.

\medskip

\setlength{\tabcolsep}{5pt}
  \noindent\begin{tabular}{|l|l ll l l l l l l l |}
    \toprule
   No & $0$ & $1$ & $ 2$ & $.....$ & $ {m-5}$ & $ {m-4}$ & $ {m-3}$ & $ {m-2}$ & $ {m-1}$ & $ m$\\
    \toprule
   $x:$ &-3 & $6q_1-3$ & $-1$ & ... & -1 & -1 & $-0.92$ & $-0.899$ & $-0.447$ & 0 \\
 \midrule
     $y:$ &0.8 & 0.8 & 0.8 & ... & 1.750 & 1.8 & 1.84 & 1.798 & 0.894 & 0 \\
 \midrule
     $q_i:$ &- & $q_1$ & $\dfrac{1}{3}-q_{1}$ & ... & ... & 0.008 &0.015 & 0.008 & 0.168 & $\dfrac{1}{6}$ \\
 \midrule
     $\al_i:$ &- & $\al_1$ & $0$ & ... & ... & 0 &0 & 0& 0 & $0$ \\
    \bottomrule
  \end{tabular}
\medskip

\setlength{\tabcolsep}{6pt}
 \noindent\begin{tabular}{|l|l ll l l l l l l  |}
    \toprule
  No &  $ {m+1}$ & $ {m+2}$ &  $ {m+3}$  & $ {m+4}$ & $ {m+5}$ &......& $ {2m-2}$ & $ {2m-1}$ & $ {2m}$\\
    \toprule
   $x:$ & $\rho\sin\te$ & $0.899$ & 0.92 & 1 & 1 & $...$ & $1$ & $3-6q_{2m}$ & 3 \\
 \midrule
     $y:$ & $\rho\cos\te $& $1.798$ & 1.84 & 1.8 & 1.75 & ... & 0.8 & 0.8 & 0.8 \\
\midrule
     $q_i:$ &$\dfrac{\rho}{6}$ & $q_{m+2}$ &  0.008 & 0.015 & 0.008 &.... & .... & $\dfrac{1}{3}-q_{2m}$ & $q_{2m}$ \\
 \midrule
     $\al_i:$ &$\fy$ & $ 0$ & 0 & $\al_{m+4} $& 0 &.... & .... & $0$ & $\al_{2m}$ \\
    \bottomrule
  \end{tabular}

\medskip
\bigskip

{\bf Some  comments  on the  values  in the  table.}

\medskip

1. Since the  most significant  point  is $z_m$, it lies  in the  origin. Therefore $z_0=(-3,0.8,0)$ and  $z_{2m}=(3,0.8,0)$. 

2. {$ \bf q_1, \al_1, q_{2m}, \al_{2m}$} :   According to  the  Proposition \ref{dense in C2},  the  pairs  of
generators of  dense subgroups  of  the  second  type are  dense
in $\bbc^2$. So we take such pair  of  generators $q_1e^{i\al_1},
q_{2m}   e^{i\al_{2m}} $ lying in   $0.003$-neighborhood of  the
point 1/6.  $q_1$ and $q_{2m}$  define  the  points $z_1$, $z_{2m}$  and  corresponding  contraction  ratios.

3. {$ \bf q_{m+1} , \al_{m+1},  q_{m+2}$}:   The point $z_{m+1}$ has  the  norm $\rho$ and  the  angle  between $z_mz_{m+1}$ and OY axis  is $\te$, this gives  us  the  values  in the column $m+1$.  $q_{m+2}$ is  defined  as 
$\dfrac{\|z_{m+2}-z_{m+1}\|}{6}$.

4. {\bf Symmetry.}   For any $i\neq 1,\  m-1,\ m+1,\ 2m-1$,  the  points $z_i$ and $z_{2m-i}$ are  symmetric  with respect  to the y axis.

\begin{figure}[h]  {\centering\includegraphics[scale=2.2]{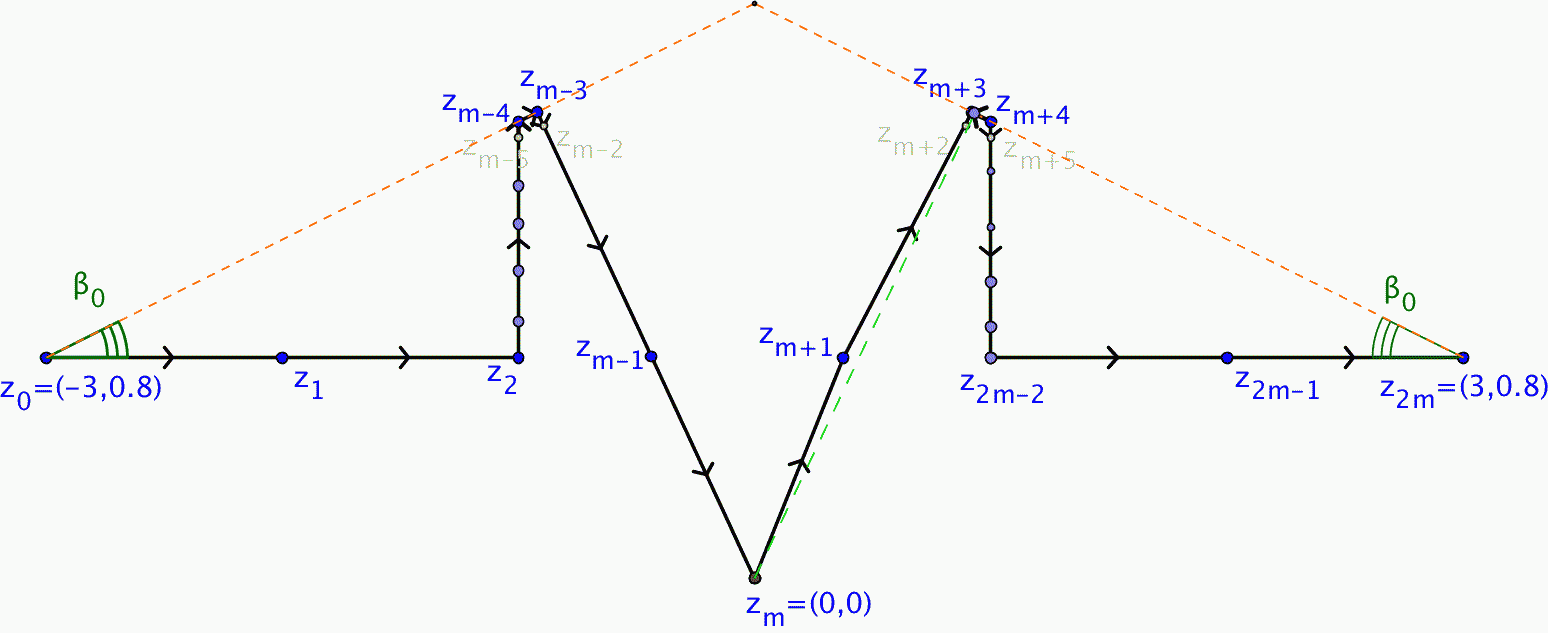}}
\end{figure}

5. {\bf $\be_0$ triangle.}    The points $z_{m-4},z_{m-3},z_{m+4},z_{m+3}$ lie  on the  sides  of  an isosceles  triangle with angles $\be_0=\arctan(1/2)$ at its base $z_0z_{2m}$. \\

\begin{figure}[h] {\centering\includegraphics[scale=2.3]{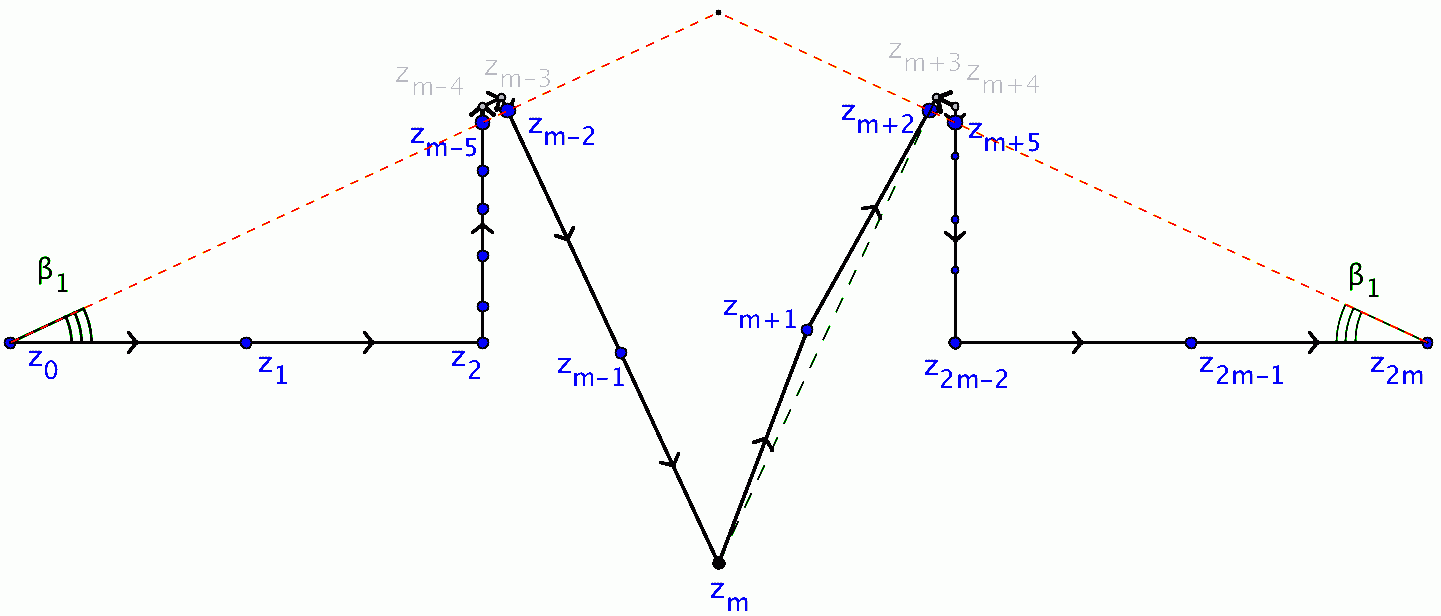}}
\caption{}
\end{figure}

  6. {\bf$\be_1$ triangle.}   At the same  time the points $z_{m-5},z_{m-2},z_{m+2},z_{m+5}$ lie  on the  sides  of  an isosceles  triangle with angles $\be_1$  at its base $z_0z_{2m}$. \\

7. {\bf Points  near $z_{m}$.}   The points $z_{m-3}, z_{m-2}, z_{m-1}, z_{m}$ lie on a line, forming the angle $\be_0$ with  the vertical axis and $||z_{m}-z_{m-1}||=1$. The same is true for their symmetrical counterparts $z_{m+3}, z_{m+2},  z_{m}$. The  point $z_{m+1}$ is  slightly  shifted to the  left and upwards.

8.{ \bf Omitted  entries.} The points $z_3 ,..., z_{m-6}$ and $z_{m+6} ,..., z_{2m-3}$ divide
the intervals $(z_2 , z_{m-5})$ and $(z_{m+5} , z_{2m-2})$ into
sufficiently  small  equal  parts, therefore  their coordinates
and  corresponding $q_i$  need  not  be  mentioned.

\subsection{The  similarities $S_i$}

\begin{figure} [h] {\centering\  \includegraphics[scale=.28]{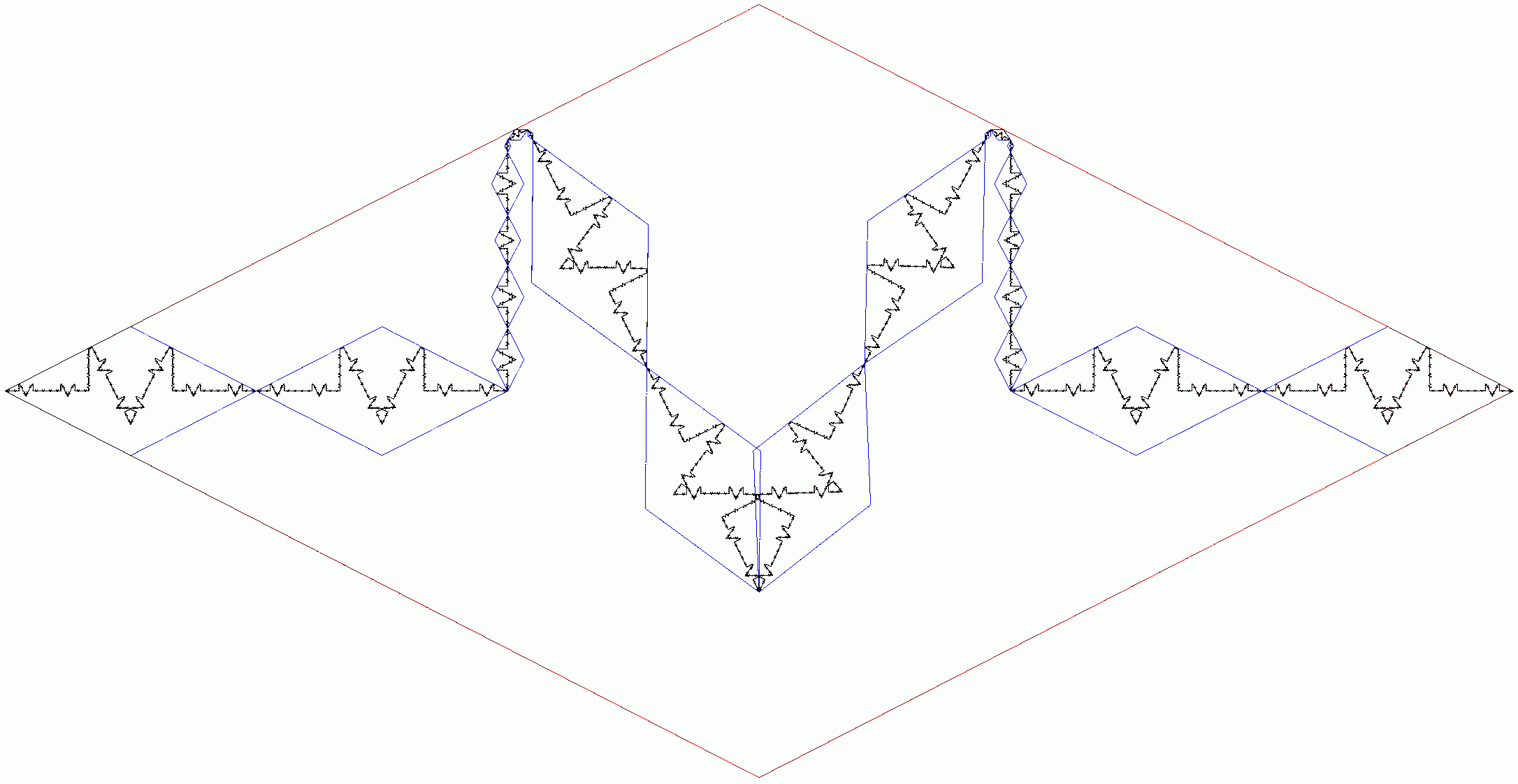}}
\caption{The  attractor $\ga$ of al  zipper $\eS^\xi$ and  the  sets $V_i$.}
\end{figure}

\begin{figure}[h]  \includegraphics[scale=.21]{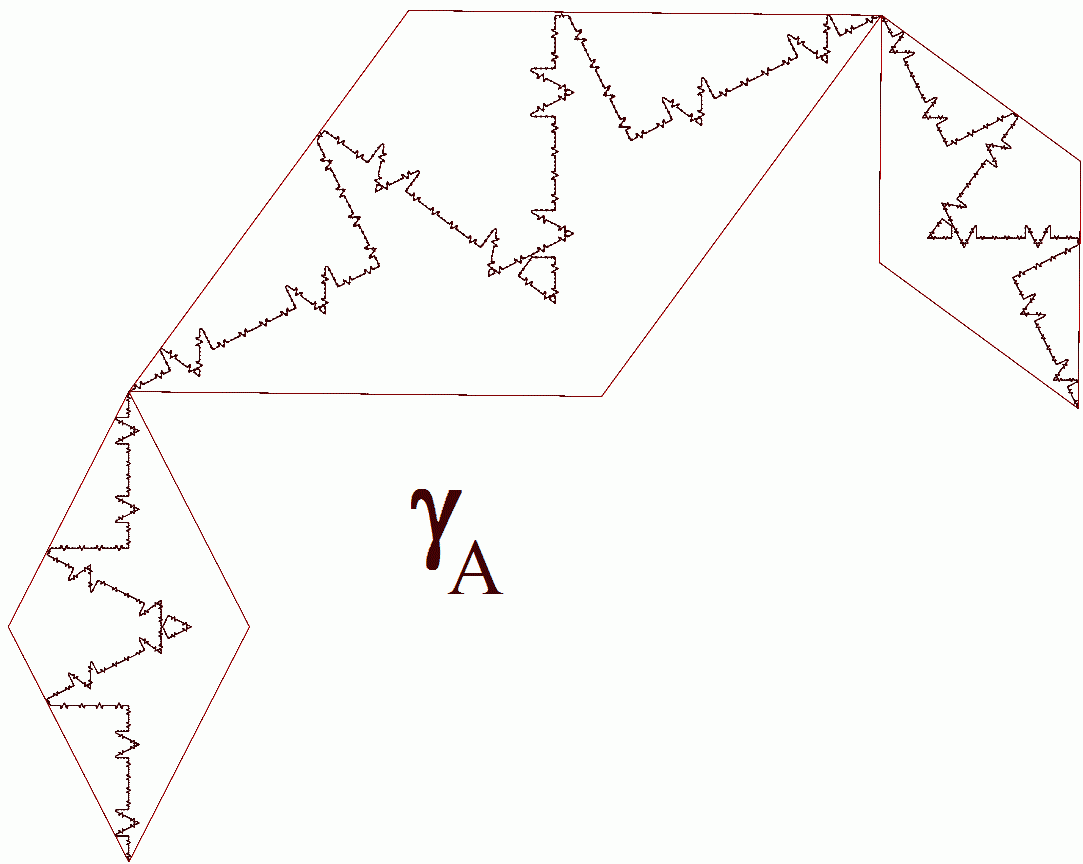}\  \  \ \    \includegraphics[scale=.21]{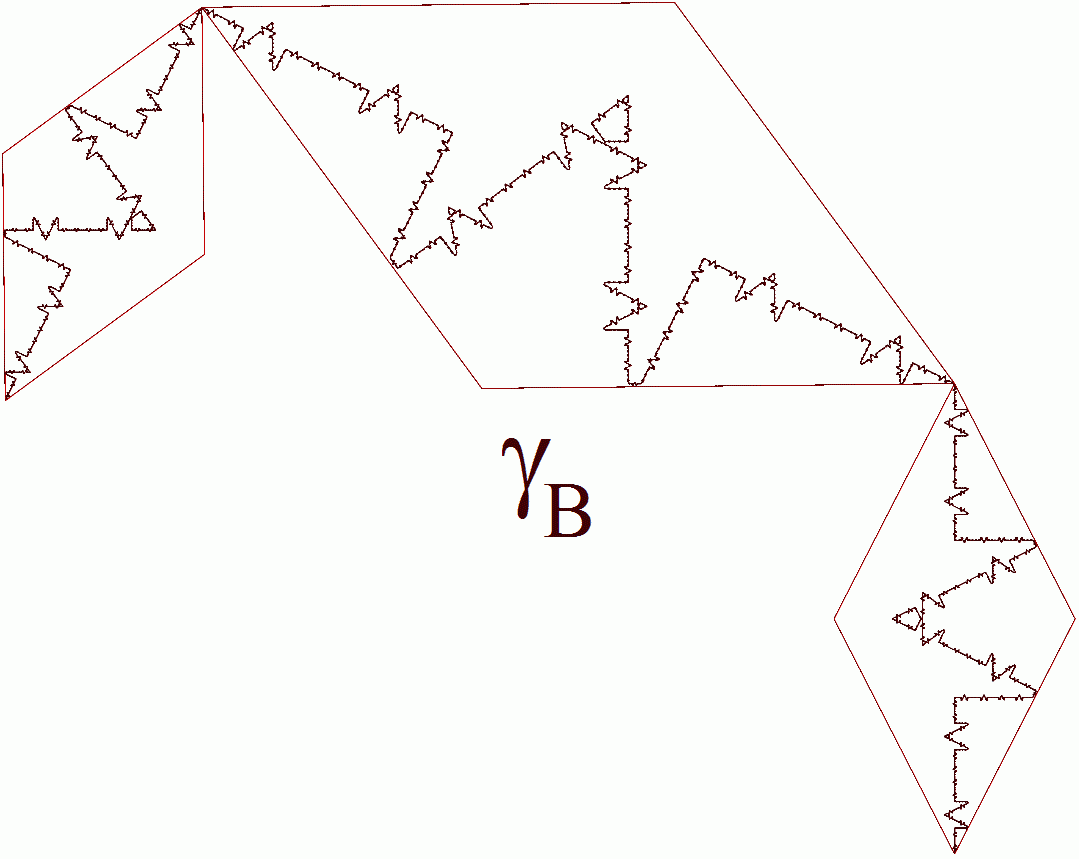}\quad \includegraphics[scale=.26]{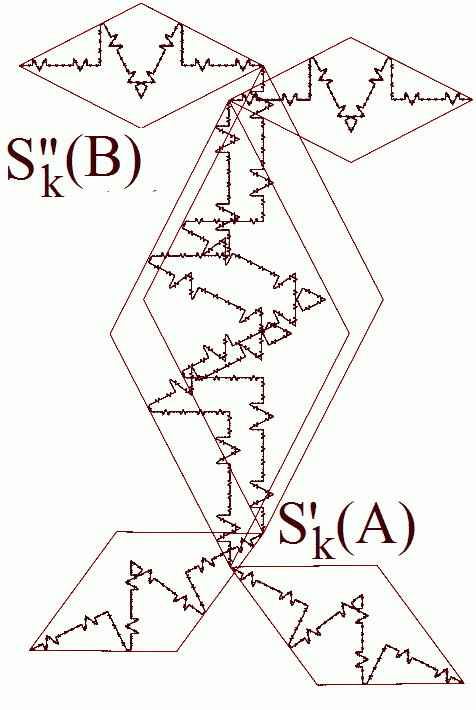}\\
\caption{The  sets $\ga_A=\ga\cap A$ and $\ga_B=\ga\cap B$. The position of $\ga_{m+4}$  is due  to  signature $\e_{m+4}=1$.  So  their images $S'_k(\ga_A)$ and $S_k"(\ga_B)$ are parallel.} \label{ABreverse}
\end{figure}

Each zipper $\eS^\xi=\{S_1,....,S_{2m}\}$  has  vertices $(z_0,...,z_{2m})$ and  signature\\ $(0,0,....0,1,0,...0)$,  where  only $\e_{m+4}=1$. This means  that  $S_{m+4}$ reverses the  order, i.e.  $S_{m+4}(z_0)=z_{m+4}$ and $S_{m+4}(z_{2m})=z_{m+3}$. See Fig. \ref{ABreverse}. 

$S_1$ and $S_{2m}$ are defined  as
 the compositions of   contractions with ratios $q_1, q_{2m}   $ with fixed  points $z_0$  and $z_{2m}   $ respectively  and  rotations  in angles $\al_1$  and   $ -\al_{2m}   $  around  the  real axis.

The map $S_{m+1}^\xi $  is a  composition of a plane similarity sending $ z_{0}$ to $ z_{m} $ and  $z_{2m}   $ to $z_{m+1}$ and  a rotation around the line $z_{m}  z_{m+1} $ in the angle $\fy$, which is the third coordinate of the parameter $\xi$.

The  map $S_{m+2}$  preserves the plane $XY$, but depends on $\xi$ because its value at $z_0$ is $z_{m+1}^{\xi}$.

Finally, the map $S_{m+4}$ is  a composition of a (fixed) similarity   map of a plane $XY$ sending $z_{0}$ to $z_{m+4}$ and $z_{2m}   $ to $z_{m+3}$ and  a rotation in an  angle $\al_{m+4}$ around the line
 $z_{m+3}z_{m+4}   $. The  angle $\al_{m+4}$ depends on  the  coordinate $\te$ of $\xi$ and will be  defined  in the  next  section.

All the other  maps $S_2,...,S_{2m-1}$ are the  similarities, preserving  XY plane.

\begin{wrapfigure}[11]{r}[40 pt]{5cm}
\includegraphics[scale=.48]{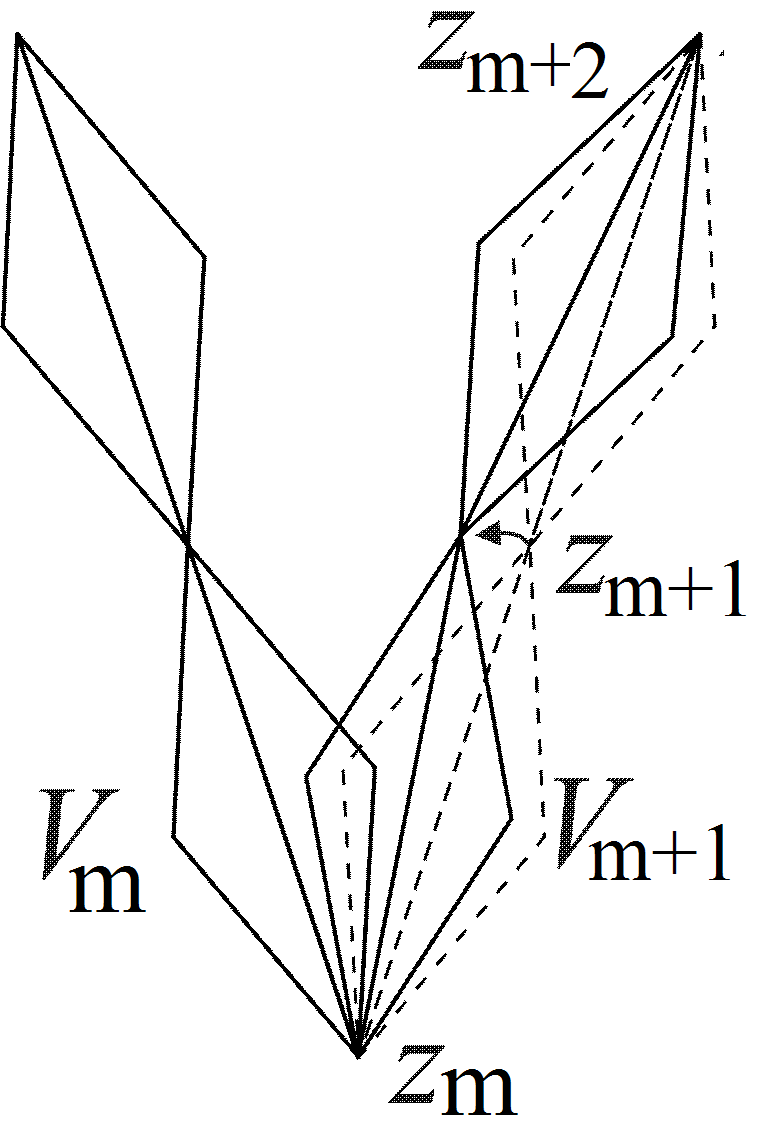}\caption{}\label{deform}
\end{wrapfigure}

The  similarity ratios $q_i$ of $S_i$ are  equal to
$|z_i-z_{i-1}|/6$. Direct  computation shows  that  when $m\ge 12$ and $q_3=...=q_{m-5}$, the  similarity
dimension  of $\eS$ is  less than 1.28 for any $\xi\in D$.

\subsection{ Bicones  and  the  sets $A$ and  $B$.}

Let $V'$ be  a  rhombus  whose  diagonal is  $z_0, z_{2m}   $ and
the  angle between diagonal  and its  sides is  $\be_2$ .

   The  value  of $\be_2$ was  chosen  as  minimal of those ones, 
 for which  two small
copies of $V'$ with  diagonals $z_{m-3}z_{m-4}$ and
$z_{m+4}z_{m+3}$  respectively lie  inside  the large $V'$. 

\begin{figure}[h] \includegraphics[scale=1.1]{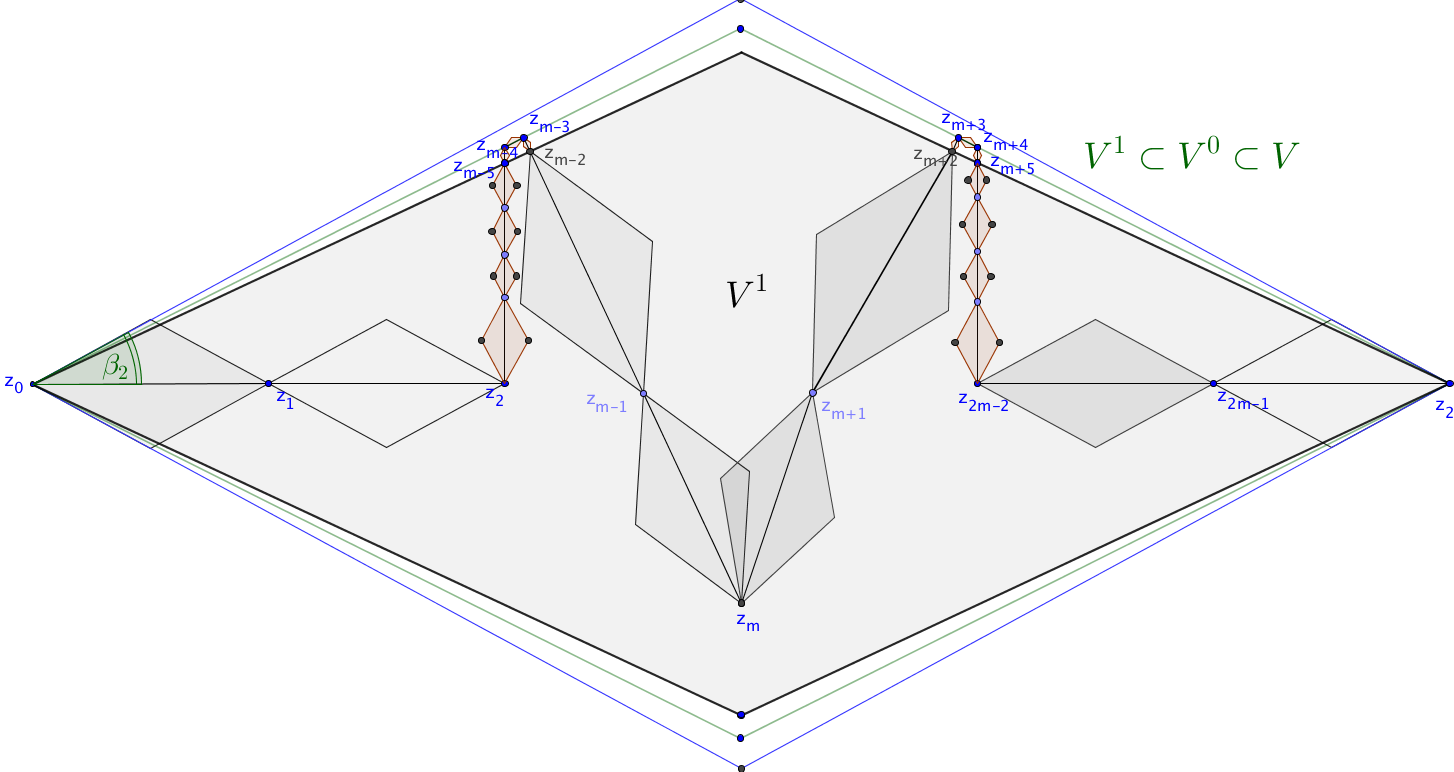}
\caption{The bicones $V^1\IN V^0\IN V$: $\be_2 $ is enlarged to make the  inclusion visible.}
\end{figure}

For  each  other $i$,  the copy of $V'$ whose  diagonal is $z_{i-1}z_i$ also  lies  inside $V'$.

Getting  to  the  space,   we  replace   $V'$ by a bicone $V$ with the  same axis $z_0,z_{2m}   $  and  angle $\be_2$ between  the  axis  and   a generator.

\medskip

 We denote by $V_i$ the  images of the  bicone $V$  under  similarities  sending
$z_0, z_{2m}   $ to $z_{i-1},z_i$.   They satisfy the  relations {\bf (A1)--(A3)} from Section 1. 

\medskip

%
%
%
%

These  properties  are  valid  for any  choice  of $\xi\in D$. (See Fig.\ref{deform})

Along with the bicones $V_i$, we  will consider the  bicones $V_i^0$ with the  same  vertices  and   with  the  angle $\be_0$ between axis  and  generator and  the  bicones $V_i^1$ with the  same  vertices  and   with  the  angle $\be_1$ between axis  and  generator.

\begin{lem}\label{allaboutV}
For  any $\xi\in D$,  \\
(i)  for  any $i=1,...,2m$, $V^0_{i-1}\cap V^0_{i }=\{z_i\}$;\\
(ii)  $V_m\cap V^1_{m+1}=V^1_m\cap V_{m+1}=\{z_m\}$;\\
(iii) the dihedral  angle  with  edge $z_mz_{m+1}$ containing common generators  of $V_{m}$ and  $V_{m+1}$, is  no greater  than $0.545$\\
(iv)   $V^0_{m}\cap V^0_{m+1}\neq\0$ and the dihedral  angle  with the  edge $z_mz_{m+1}$ and  sides, containing common generators  of $V^0_{m}$ and  $V^0_{m+1}$, lies  between $0.224$ and $0.317$.

\end{lem}

\dok   \ \ Since  the  angle between  the  axes  is  greater  than  $2\be-\mu$, we  have (i) and (ii).\\
The  upper  bound   of  such  angle  for $V_{m}\cap V_{m+1}$ is $2 \arccos\dfrac{\tan(\be-\mu/2)}{\tan(\be+\mu)}=0.545$ gives  us (iii). Now, the  upper  and  lower  bounds  for the  angle  for 
$V^0_{m}\cap V^0_{m+1}$ are  $2 \arccos\dfrac{\tan(\be-\mu/2)}{\tan(\be)}$ and $2 \arccos\dfrac{\tan(\be-\mu/4)}{\tan(\be)}$, which  gives (iv).\vse

\bigskip

\subsection{ The  value  of $\al_{m+4}$ .}\label{amplus4}

Consider the bicones $V^0_{m}  $  and $V^{0}_{m+1}$. They  have  common  vertex $z_m=0$. The  bicone $V^0_m$ stays fixed, while $V^0_{m+1}$ changes  its position  depending on its  second  variable  vertex $z^\xi_{m+1}$, and  the angle  between  the  axis of $V_{m+1}$ and OY is $\te\in(\be-\mu,\be-\mu/2)$, therefore $V^0_{m}  $  and $V^{0}_{m+1}$ have  nonempty  intersection. 

\begin{wrapfigure}[12]{r}[30 pt]{5cm}
\includegraphics[scale=.5]{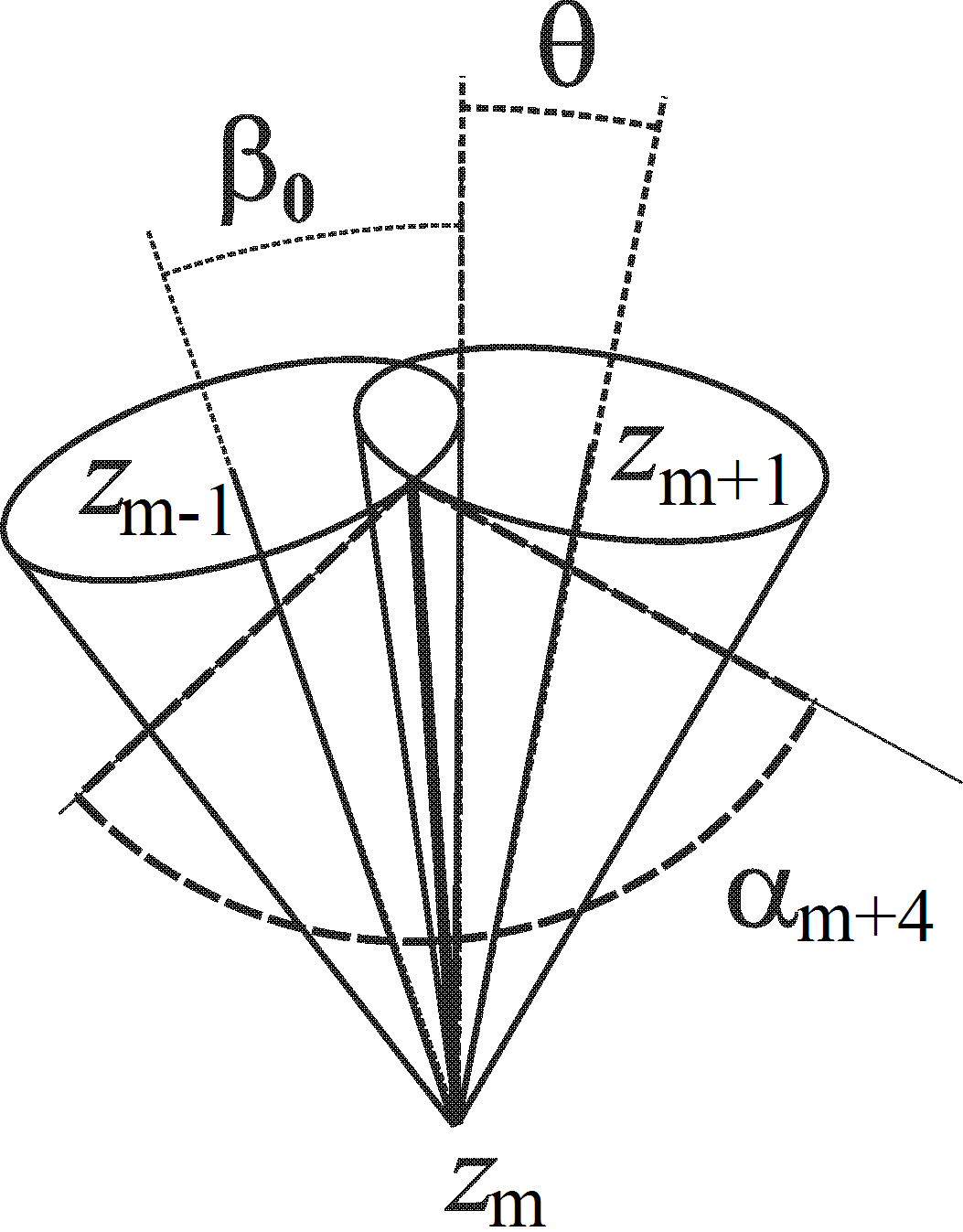}\label{angle}
\end{wrapfigure}  The angle $\al_{m+4} $ is equal to an angle  between the  tangent  planes   to  the  surfaces of $V^0_{m}  $  and $V^0_{m+1}$ at    points of  their  common generator. By spherical cosine  theorem, it
is  equal to $$\al_{m+4}(\te)=\arccos(4-5\cos(\be_0+\theta)). $$

\bigskip

The values of the derivative  $\cos(\al_{m+4}(\te))'$   for $\theta\in[\be-\mu,\be-\mu/2]$ lie  in  the  interval $(-20,-14)$. Therefore, for  any $\te_1, \te_2 \in (\be-\mu,\be-\mu/2)$,  $$|\al_{m+4}(\te_1)-\al_{m+4}(\te_2)|<20|  \Da\te  |$$

\subsection{ The  sets A  and B.}

\bigskip
  
\begin{figure}[h] \quad \quad  \includegraphics[scale=.26]{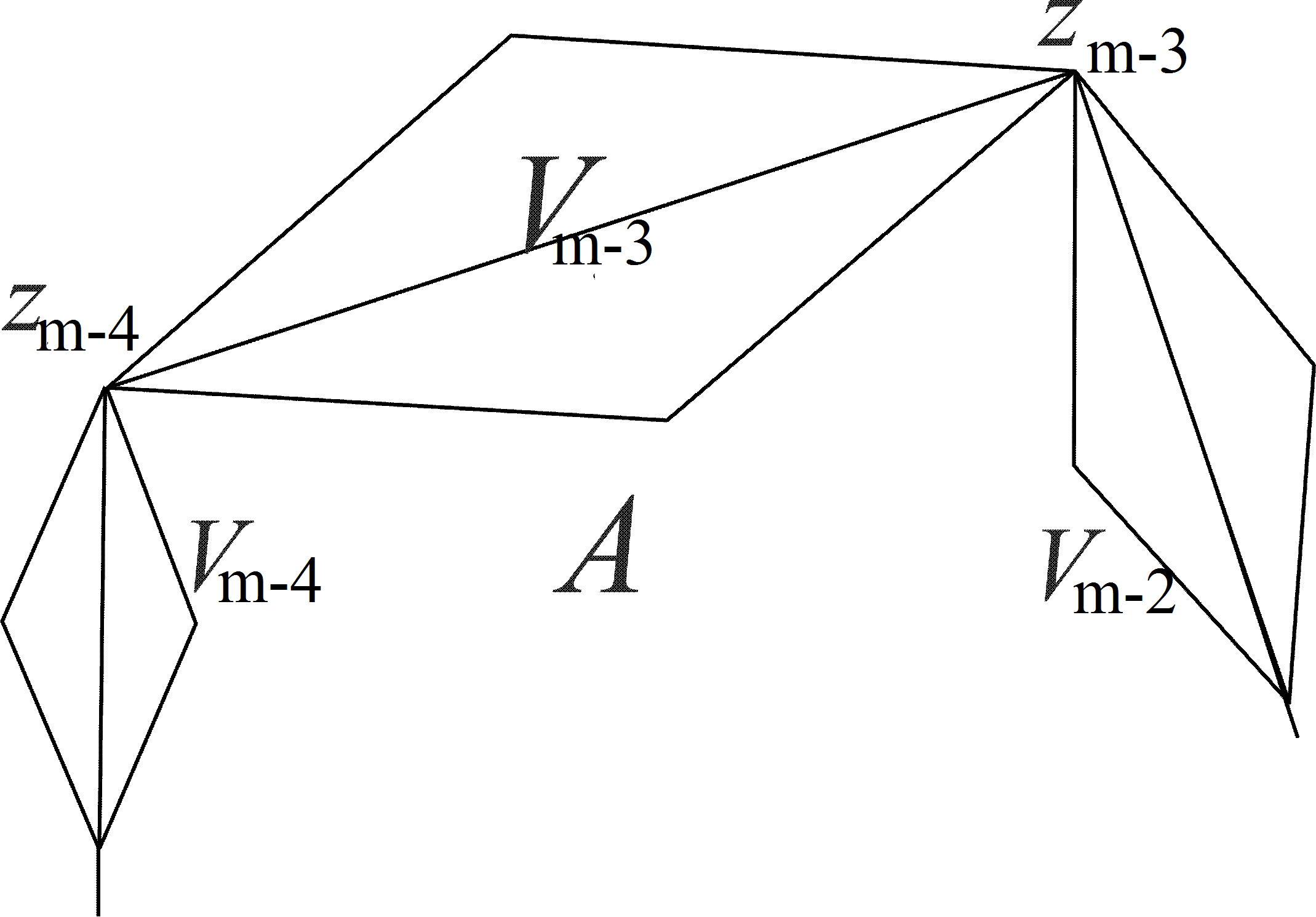}\ \quad \quad  \quad \quad \ \    \includegraphics[scale=.26]{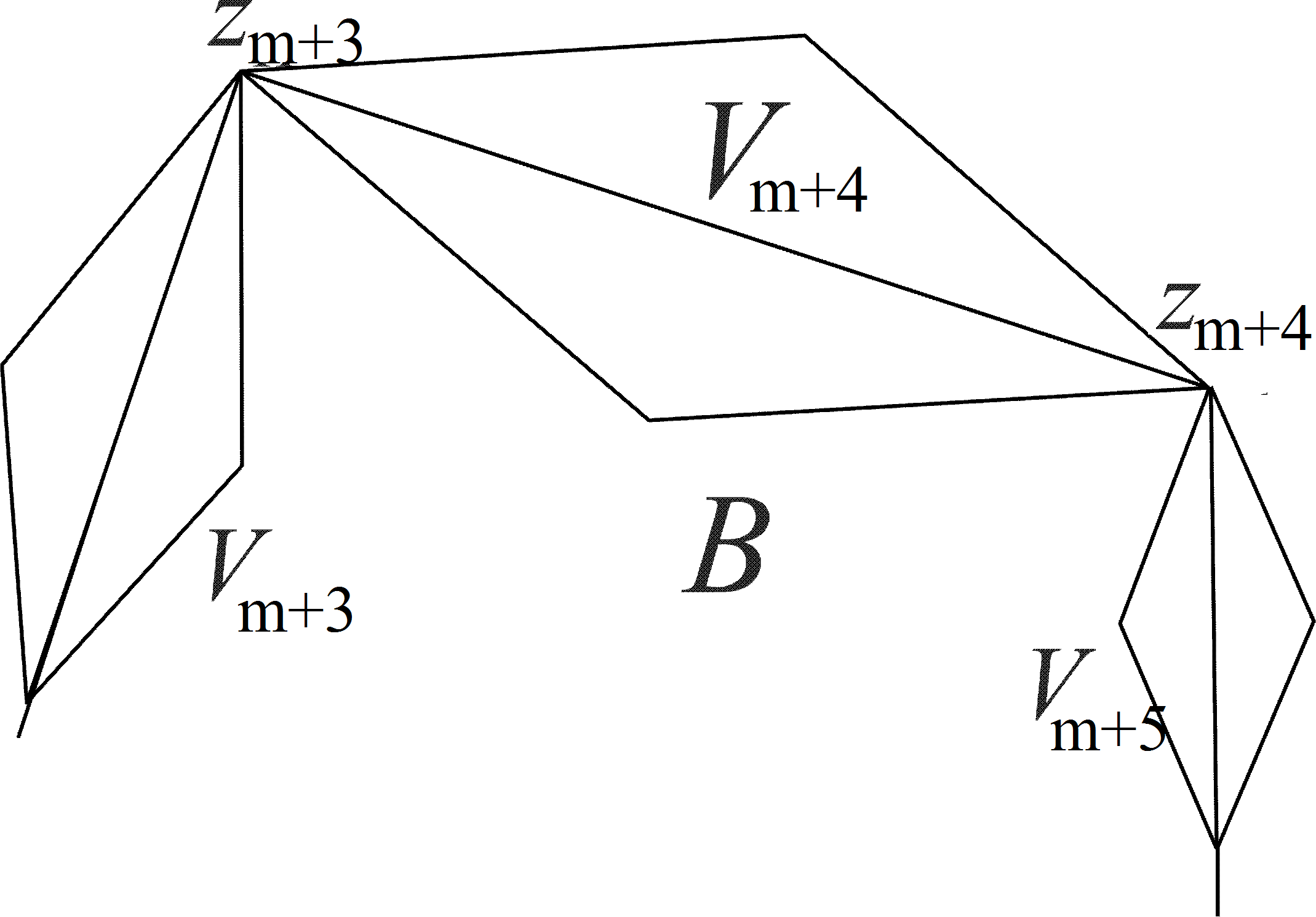}\\
\caption{The  sets  A and B.}
\end{figure}

 There  are 2 sets, formed  from $V_i$-s, which will be needed  for  our  considerations: the  set  $A=V_{m-4}\cup V_{m-3 }\cup V_{m-2}$ and  the  set
$B=V_{m+3}\cup V_{m+4 }\cup V_{m+5}$:

\bigskip

The  sets $A$ and $B$ lie outside the bicone $V^1$, but intersect $V^0$ so that  the  axes  of $V_{m-3}$ and $V_{m+4}$ are contained in the boundary of $V^0$ .

 { \bf  Some  inequalities  for  the  sets  $A$ and $B$.  }

For any  two points $x$ and $y$ in  the  set $A$ (or  both in the  set  $B$),  we have $\dfrac{||x-z_0||}{||y-z_0||}<1.06$ and  $\dfrac{||x-Pr_1(x)||}{||y-Pr_1(y)||}<1.095$. At the  same time, the  dihedral  angle, containing A and B, whose edge is  the line $z_0z_{2m}$,  is  no  greater  than $2\sqrt{5}\mu$.

Therefore,  the  set A lies  in  a domain defined by  the  inequalities
$$ R\le |x-z_0|\le 1.06 R;   \ \ \ \    -\sqrt{5}\mu\le \fy\le\sqrt{5}\mu; \ \ \ \    \be_0-2\mu\le \te\le\be_0+\mu,$$
where $ \fy$ and $ \te$ are the  azimuth  and polar  angles for  the
point $x$ in  spherical  coordinates  with the origin $z_0$, real
line being  the polar  axis and the  azimuth  direction being
parallel to  $OY$.

Here $R=2.214$ is  the  distance   from $z_0$ to the  nearest  point  of  $A$.\\

Direct  computation  shows that

\begin{lem}  The set $A$  can be covered by  a ball  W of radius $0.036R$ with the  center  defined by

$$  |x-z_0|= 1.03 R;   \ \ \ \     \fy=0; \ \ \ \     \te=\be_0-\mu/2 \quad \blacksquare$$   \end{lem}

By symmetry,  analogous inequalities hold for the points in the set $B$.\\

\section{ Checking  the  properties  of  the  family  of  zippers $\eS_\xi$.}

Now  we  can  verify  the  properties $\bf A1-A6$ for  any $\xi\in D$ and evaluate the  difference $S_i^\xi(x)-S_i^\eta(x)$ for $x\in V$. Fortunately {\bf A1-A3} are  obvious so we  proceed to {\bf A6, A5, A4}.

\subsection { Defining  linear  zipper $\eT$ and  checking A6.}

As it was proved in ([1], Lemma 1.1), for  any  linear  zipper $\eT=\{T_1,...,T_{2m}\}$ with attractor $K(\eT)=[0,1]$ and with  the same signature $\e$ as  the  zipper $\eS$, there  is  unique   continuous function 
$\fy_{\eS\eT}:[0,1]\to \ga(\eS)$ satisfying 
$$ \fy_{\eS\eT}\circ T_i(t)=S_i\circ \fy_{\eS\eT}(t) \mbox{ for  any  } i =1,...,2m\mbox{  and  } t\in[0,1]$$
This  map is  H\"older  continuous. We prove  the  following theorem which  allows  to  find its  H\"older  exponent:

\begin{thm}\label{Holder}

Let  $\eS=\{S_1, \dots, S_m\}$  be  a  zipper  in $\mathbb{R}^n$ and  suppose all  $S_i$ are   similarities.\\
Let $\fy$: $I[0,1]\rightarrow \ga(\eS)$ be the linear parametrization of $\eS$ by a zipper $\eT=\{T_1,\dots T_m\}$. The H\"older  exponent $\al$ of  the map $\fy$ satisfies $$\alpha = \min \limits_{i=1,\dots, m} \dfrac{\log \Lip  S_i}{\log\Lip T_i}$$

\end{thm}

\dok

Denote $q_{i} = \Lip  S_{i}$, $p_{i} = \Lip  T_{i}$, $p_{\min}=\min \limits_{i=1,...,m} p_i$. 
Let  $M$ be  the  diameter  of $\ga(\eS)$.

Observe  that  for  any  multiindex $\bi={i_1...i_k}$, 
$q_\bi\le p_\bi^\al$.
\probel

Suppose $a, b \in I$, $a < b$. There  are  two  possibilities:

1) For  some  multiindex $i_1...i_ki_{k+1}$, $$T_{i_1...i_ki_{k+1}}(I)\IN [a,b]\IN T_{i_1...i_k}(I)$$
In this  case, $p_{i_1 \dots i_k} p_{\min} < |a-b| \leq p_{i_1 \dots i_k}$, while  $\|f(a)-f(b)\| \leq q_{i_1 \dots i_k}M$. 

But  $ q_{i_1 \dots i_k}\le  p_{i_1 \dots i_k}^{\al}$, so   $$\|f(a)-f(b)\| \le {\dfrac{M}{p_{\min}^{\al}}} {|a-b|}^{\al}$$

2) For  some  pair of multiindices $i_1...i_ki_{k+1}$ and $j_1...j_lj_{l+1}$, 
$$T_{i_1...i_ki_{k+1}}(I)\cup T_{j_1...j_lj_{l+1}}(I) \IN [a,b]\IN T_{i_1...i_k}(I)\cup T_{j_1...j_l}(I) $$

In this  case, $(p_{i_1 \dots i_k} +p_{j_1...j_l})p_{\min} \le |a-b| \leq (p_{i_1 \dots i_k} +p_{j_1...j_l})$, while \\ $\|f(a)-f(b)\| \le (q_{i_1 \dots i_k}+q_{j_1...j_l})M$. 

Suppose  $p_{i_1 \dots i_k} \ge p_{j_1...j_l}$, then   $p_{i_1 \dots i_k} p_{\min} \le |a-b| \leq 2p_{i_1 \dots i_k} $

Thus, we  have  $$\|f(a)-f(b)\| \le 2M p_{i_1 \dots i_k}^\al \le \dfrac{2M}{p_{\min}^\al}|b-a|^\al$$

The  exponent $\al$ is  exact  because for  some k,  $q_k=p_k^\al$, and  therefore 
$\|f(a)-f(b)\|= L {|a-b|}^{\al}$ for $a=T^n(0), b=T^n(1), L=\|z_0-z_m\|$. \vse

\begin{lem}\label{oneforall}
There  exist  such linear  zipper $\eT$ on $[0,1]$, that  for  any $\xi\in D$, the  linear  parametrisation $\fy^\xi=\fy_{\eS^\xi\eT}$  has  H\"older  exponent greater  than 3/4.

\end{lem}

\dok \ \ 
Take some $ \xi_0=(1,0,\te_0 )\in D$. Let $q_i=\Lip S_i^{\xi_0}$ and $s$ be  the  similarity  dimension of $\eS^{\xi_0}$. 
Take  a zipper $\eT=\{T_1,...,T_{2m}\}$ on  [0,1] with  signature $\e$ and  contraction  ratios $p_i=q_i^s$. Obviously $\dfrac{\log q_i}{\log p_i}\equiv 1/s $. There  are  only  two indices $i=m+1$  and $i=m+2$, for  which $\Lip S^\xi_i=q'_i$ depends  on $\xi$. For  both of  them

  $|\log q'_i/q_i|\le \log(1.02)$, therefore $1/1.012<\dfrac{\log q'_i}{\log q_i}<1.012$.

Therefore $  \dfrac{\log p_i}{\log q'_i}< 1.28\cdot 1.012<4/3$.

It follows  that  for  any  $\xi\in D$, the  linear parametrization $\fy^\xi$of $\eS^\xi$ by  the  zipper $\eT$ has  H\"older  exponent greater  than $3/4$. \vse

\subsection{  Verification of A5}

First, we reformulate the  first statement  of Theorem 4 in the following way:
\begin{lem}\label{tworays}
 Let $\eta_1, \eta_2$  be the generators  of  a dense subgroup  of  second type, then for  any   $z_1,z_2\in \mathbb{C}\setminus\{0\}$ and  any rays $l_1,l_2$ issuing from  zero, there is such sequence $\{(n_k,m_k)\}$  that $\lim\limits_{k\to \infty}   \left|\dfrac{z_1\eta_1^{n_k}}{z_2\eta_2^{m_k}}\right|=1$,  while the angles between $l_1, l_2$ and rays  passing  from  0 to ${z_1\eta_1^{n_k}}$  and ${z_2\eta_2^{m_k}}$ respectively, converge  to 0. \vse \\
\end{lem}

Second, the statement of  the  Lemma  remains  valid if we replace:\\
 1) the plane $\bbc$ by a  cone $C\IN\rr^3$ with  the axis $L$; \\
2) the products $z_1\eta_1^{n_k}$ and $z_2\eta_2^{m_k}$ by
$f_1^{n_k}(z_1)$ and $f_2^{m_k}(z_2)$, where $f_i$  are compositions of homothety with contraction  ratio $|\eta_i|$ and rotation around the axis $L$ axis in  the  angle $arg(\eta_i)$   and\\
 3) Taking  for  $l_1,l_2$ some  generators of  the  cone $C$.

\medskip

Since  for any two  sequences $f_k$, $g_k$ of  orthogonal transformations in $\rr^n$ (and  therefore for similarities  having origin as a  fixed  point), the  convergence  of the  sequence
$f_k^{-1}g_k$ to ${\rm Id}$ is  equivalent to the  convergence of $g_k f_k^{-1}$ to ${\rm Id}$, we  have the  following

\begin{lem}\label{segmaps} Let $PQ$ be  a  segment in $\rr^3$ and  $\bn$ be some its  normal vector.
Let $D_k, \ D'_k$ be  two sequences of  points  in $\rr^3\mmm \{0\}$ and $\bn_k,\bn'_k$ sequences  of  unit normal  vectors  to the segments $OD_k,\  OD'_k$. Let $f_k$ (resp. $g_k$)
be  the  similarities which map  $P$ to $O$, $PQ$ to $OD_k$  (resp. $OD'_k$) and  whose orthogonal parts send $\bn$ to $ \bn_k$ (resp. $\bn'_k$). Then \\

\quad\quad$f_k^{-1}g_k \to {\rm Id}$\quad\quad iff \quad\quad $\dfrac{|OD_k|}{|OD'_k|}\to 1$,  $(\bn_k,\bn_k')\to 1$, $\angle D_k OD'_k\to 0$. \vse
\end{lem}

\begin{lem}
For  any $\xi\in D$,  there  is  such   sequence $\Sa_\xi=\{(i_k,j_k)\}$ that the  similarities
$S'_k=S_{m}S_{2m}^{i_k}$ and $S''_k=S_{m+1}S_{1}^{j_k}$   satisfy
$$\lim\limits_{k\to\infty}(S'_{k}S_{m+4})^{-1}(S''_{k}S_{m-3})={\rm Id}.$$
\end{lem}

\dok     \ \ \  Observe  that $S_{m+4}^{-1}(z_{2m})=S_{m-3}^{-1}(z_0)=(-153,0.8,0)$. Denote  this  point  by $P$.
Let $ Q=(3,0.8,0)$. Each map $S_{m}S_{2m}^{k}S_{m+4}$ sends $P$ to the point $z_m=(0,0,0) $ and the point $Q$ to  some point on the surface of the cone $S_m(V^0)$. Each map $S_{m+1}S_{1}^{k}S_{m-3}$ sends $P$ to the point $z_m=(0,0,0) $ and the point $Q$ to  some point on the surface of the cone $S_{m+1}(V^0)$. The maps $f_1=S_mS_{2m}S_m^{-1}$ and $f_2=S_{m+1}S_{1}S_{m+1}^{-1}$ are the similarities preserving the  cones $S_m(V^0)$ and $S_{m+1}(V^0)$ respectively, defined by  two  generators  of  a  dense  group of  second  kind.

Let $l$ be a  common  generator  of  these  intersecting  cones.

Denoting $w_1=S_m\circ S_{m+4}(Q)$ and $w_2=S_{m+1}\circ S_{m-3}(Q)$, rewrite $$S_{m}S_{2m}^{k}S_{m+4}(Q)=f_1^k(S_m\circ S_{m+4}(Q))=f_1^k(w_1)$$ and
$$S_{m+1}S_{1}^{k}S_{m-3}(Q)=f_2^k(S_{m+1}\circ S_{m-3}(Q))=f_2^k(w_2)$$ Since $w_1, w_2$ lie  on  the surfaces of  the  cones $S_m(V^0)$ and $S_{m+1}(V^0)$, we may  apply  Lemma\ \ref{tworays}  to  get  subsequences $i_k, j_k$ for  which
 $\lim\limits_{k\to \infty}   \dfrac{||f_1^{i_k}(w_1)||}{||f_2^{j_k}(w_2)||}=1$,  while the angles between $l$ and rays  passing  from  0 to $f_1^{i_k}(w_1)$  and $f_2^{j_k}(w_2)$ respectively, converge  to 0.\\
At the  same  time, since $S_{m+4}$ contains  a  rotation in  the  angle $-\al_{m+4}$, which  is  the  angle  between normals  to  these cones  at  points  of $l$, the angle  between  the  images  of  the  normal $\bn$ to the segment $PQ$  under the maps
$S'_kS_{m+4}=S_{m}S_{2m}^{i_k}S_{m+4}$ and $S''_kS_{m-3}=S_{m+1}S_{1}^{j_k}S_{m-3}$ converges  to 0.
Therefore, by Lemma\ \ref{segmaps}, $$\lim\limits_{k\to\infty}(S'_{k}S_{m+4})^{-1}(S''_{k}S_{m-3})={\rm Id}.\quad\quad\quad\blacksquare$$

\subsection{Dividing $S_m(\ga)\cap S_{m+1}(\ga)$  to a  sequence of  disjoint pieces:  checking (A4). }

\begin{lem}\label{partition}
There  is such sequence  $\Sigma=\{(i_k,j_k)\}$ in $\nn\times\nn$ that:\\
(1) For  any $\xi\in D$, $$S_{m}  (\ga)\cap S_{m+1} (\ga)= \{z_{m}  \} \cup\left( \bigcup\limits_{i,j=1}^\infty   \left(S_{m+1} S_1^{i_k}(\ga_A)\cap S_{m}  S_{2m}   ^{j_k}(\ga_B)\right)\right)$$
(2) For  any $k$  there is such $\xi\in D$ that $$ S_{m+1} S_1^{i_k}(A)\cap S_{m}  S_{2m}   ^{j_k}(B)\neq\0$$
(3) If a  pair $(i,j)\notin\Sigma$, then for any $\xi\in D$,
$$ S_{m+1} S_1^{i}(A)\cap S_{m}  S_{2m}   ^{j}(B)=\0$$
(4) The  sequences  $\{(i_k)\}$ and $\{(j_k)\}$ are strictly
increasing  and  both projections $Pr_1:(i,j)\to i$,
$Pr_2:(i,j)\to j$ are injective  on $\Sigma$,\\
(5)  For  any $\xi\in D$, $\Sa_\xi\IN\Sa$.
\end{lem}

\dok

For  any  system $\eS_\xi, \xi\in D$, its invariant set $\ga^\xi$   is  a  Jordan  arc   if  $S_{m}  (\ga)\cap S_{m+1} (\ga)=\{z_{m}  \}$.

Consider the  set $S_{m}  (\ga)\cap S_{m+1} (\ga)$.  It is  contained  in $S_{m}  (V)\cap S_{m+1} (V)$. By Lemma \ref{allaboutV},     $S_{m}  (V^1)\cap S_{m+1} (V)=S_{m}  (V)\cap S_{m+1} (V^1)=\{z_{m}  \}$. Therefore, $$S_{m}  (V)\cap S_{m+1} (V)= S_{m}  (V\mmm \dot V^1)\cap S_{m+1} (V\mmm \dot V^1)$$

The intersection of $\ga\mmm\{z_0, z_{2m}\}$ and $V\mmm \dot V^1$ lies in the set $$\left(\bigcup\limits_{n=0}^\infty S_1^n(A)\right)\ \ \ \bigcup\ \ \  \left(\bigcup\limits_{n=0}^\infty S_{2m}^n(B)\right)$$
therefore
$$S_{m}  (\ga)\cap S_{m+1} (\ga)\IN \{z_{m}  \} \cup \left(\bigcup\limits_{i,j=1}^\infty   \left(S_{m+1} S_1^i(A)\cap S_{m}  S_{2m}   ^j(B)\right)\right)$$

Suppose $x\in A$, $y\in B$ and  $S_{m+1} S_1^i(x)= S_{m}  S_{2m}   ^j(y)$. Denote  $a=\log(\|x-z_0\|)$, $b=\log(\|y-z_{2m}   \|$.  Then $\log(q_{m+1} )+i \log(q_1)+a=\log(q_{m})+j \log(q_{2m}   )+b$, or
$$|i \log(q_1)-j \log(q_{2m}   )|\le|b-a|+\log(q_{m+1} /q_{m}  ) $$

According  to the  inequality (1), $|b-a|<\log 1.06<0.06$. At the  same time, $|\log q_{m+1}/q_m|<0.04$. Therefore $i$ and $j$ are the solutions of the inequality
$$|i \log(q_1)-j \log(q_{2m}   )|<0.1 \ \ \ \ \ \ \eqno{(**)}$$

Both $q_1$ and $q_{2m}$  lie between 1/5 and 1/7, so absolute values of their logarithms are greater than 1.

This implies that for  any $i\in\nn$ there  is  at most one    $j\in\nn$  for which the  inequality $(**)$ holds and  vice  versa. The solutions of this inequality are the same for
all values of $q_{m+1}$ satisfying    $0.98<\dfrac{q_{m+1}}{q_m}<1.02$ .

Therefore there  is a subsequence  $\{(i_k,j_k)\}$ of  $\Sa$  which  runs  through all non-negative   solutions  of  the  inequality $(**)$ for  any value  of $q_{m+1} $;  both  sequences $i_k$  and $j_k$ being  strictly  increasing. Obviously, it  will contain  $\Sa_\xi$ for  any $\xi\in D$.
 Discarding  those entries $\{(i_k,j_k)\}$,  for  which  $ S_{m+1} S_1^{i_k}(A)\cap S_{m}  S_{2m}   ^{j_k}(B)=\0$ for  any $\xi\in D$, we get   the desired sequence.

Thus,  for  any value  of $q_{m+1} $, the  set  $S_{m}  (\ga)\cap S_{m+1} (\ga)\mmm \{z_{m}  \}$is the disjoint union $\bigcup\limits_{i=0}^\infty \left( S_{m+1} S_1^{i_k}(\ga_A)\cap S_{m}  S_{2m}   ^{j_k}(\ga_B)\right)$ \vse

\section{ Proof  of  the  main  Theorem}

To prove  Theorem \ref{main} we  make  the  following  steps: 

First we make the  estimates  for  the  difference $\|S_i^\xi(x)-S_i^\eta(x)\|$ on the  set $ V$ for  any given  $\xi,\eta \in D$ for  $i=1,..., 2m$. They  are  quite  simple  for $i\neq m+1$ ( Lemma\ref{others}), so  most  of  the  work is  done  for $i=m+1$.
For each $k\in\nn$   we express  these estimates  in terms  of  the  displacement $\da^*_k=x^\xi_k-x^\eta_k$  of  the  point $x_k=S'_k(z_{m-4})$.  We prove  in Lemma \ref{Sbilip}   that  the  dependence of  $x^\xi_k$ on  $\xi$ is  bi-Lipshitz. In Lemma \ref{delta} we  estimate $\max\limits_{i=1,...,2m,x\in V}\|S_i^\xi(x)-S_i^\eta(x)\|$ in  terms  of $\da^*_k$. In Lemma \ref{bilipm+1}  we  get upper  and  lower  bounds  for  $\| S^\xi_{m+1}(x)- S^\eta_{m+1}(x)\|$ on $S_1^{i_k}(A)$ in terms  of $\da^*_k$. In the  Subsection \ref{finish}    we prove Theorem \ref{genpos} and  Proposition \ref{collage}; in Lemma \ref{fxist}  we  prove  that  the  function $f(\xi,s,t)$ used  in the  Theorem \ref{genpos}  is  bi-Lipschitz  with respect  to $\xi$  and $\al$-H\"older with respect  to 
$s,t$ which gives  us the proof  of  the main Theorem.

\bigskip

{\bf Some  notation.} \ \ \  Let $\xi=(\rho_1,\fy_1,\te_1)$, \ \ $\eta=(\rho_2, \fy_2, \te_2)$\ \  and \ \  $\Da\xi=$\\ $=(\Da \rho,\Da \fy,\Da \te)=\xi-\eta$. 
Denote  by $x^\xi_k$ the point $S'_k(z_{m-4})=S^\xi_{m+1}S_1^{i_k}(z_{m-4})$ and
let $\da^*_k=\|x_k^\xi-x_k^\eta\|$.

\bigskip

Fixing  the  value  of $\xi$, we consider  spherical  coordinate  system whose  polar  axis  is $z_mz_{m+1}^\xi$ and  whose  azimuth  direction is a  perpendicular to $z_mz_{m+1}^\xi$ lying  in the right  half-plane  of  the  plane XY. We'll denote  these  coordinates  by $\vro,\phi,\vte$ and  call them $\xi$-coordinates.

\bigskip

It is  more  convenient  to  represent    the difference $\|S_i^\xi(x)-S_i^\eta(x)\|$, $x\in A$ 
as $\|S_{m+1}^\eta(S_{m+1}^\xi)^{-1}(x)-x\|$, where $x\in S_{m+1}^\xi S_1^{i_k}(A)$. Denote $F_{\xi\eta}=S_{m+1}^\eta(S_{m+1}^\xi)^{-1}$.

Each of  the  maps  ${F}_{\xi\eta}$  may be  viewed  as  a composition  of  a  rotation $R_\fy$ in an   angle $\Da\fy$ with  respect  to  the  line $z_m z_{m+1}^\xi$,  a  rotation $R_\te$ in  an  angle $\Da\te$ with  respect to $Z$ axis and homothety $H_\rho$  with  ratio $\rho_2/\rho_1$.

\subsection{ Transition maps from $\eS^\xi$  to $\eS^\eta$  and  their  estimates.}

 First  we consider all $i\neq m+1$:

\begin{lem}\label{others} For  any $x\in V$, $\xi,\eta\in D$,\\ 
(a) If $i\neq m+1, m+2$ or $m+4$, then $S_i^\xi(x)-S_i^\eta(x)\equiv 0$;\\
(b) $\|S_{m+4}^\xi(x)-S_{m+4}^\eta(x)\|<0.46|\Da\te|$;\\
(c)  $\|S_{m+2}^\xi(x)-S_{m+2}^\eta(x)\|\le \|z_{m+1}^\xi-z_{m+1}^\eta\|$

\end{lem}

\dok.  \ \ \  For $i=m+4$, $$||S_{m+4}^\xi(x)-S_{m+4}^\eta(x)||\le q_{m+4}|\al^\xi_{m+4}-\al^\eta_{m+4}|r_x, $$ where $r_x$ is  a distance  from  the point  $x$ to  the  OX  axis. 
Therefore,  $$||S_{m+4}^\xi(x)-S_{m+4}^\eta(x)||<20 q_{m+4}|  \Da\te  |r_x<0.298|  \Da\te  |r_x$$
Thus, $r_x\le 1.538$  implies   $\max\limits_{x\in V}||S_{m+4}^\xi(x)-S_{m+4}^\eta(x)||<0.46|  \Da\te  |$.
\bigskip

The  similarity $S_{m+2}^\xi$ depends  only  on  a position  of  the  point $z_{m+1}$, therefore $$\|S_{m+2}^\xi(x)-S_{m+2}^\eta(x)\|\le \dfrac{\|x-z_{2m}\|}{\|z_0-z_{2m}\|}\|z_{m+1}^\xi-z_{m+1}^\eta\|\quad \quad \blacksquare$$

\bigskip

The  estimates  for $i={m+1} $  are  more  complicated.

\bigskip

The  following  two lemmas   will help  us  to  estimate  maximal  and  minimal  displacement  for different  points  in $V_{m+1}$:

\begin{lem}\label{inaball}
Let $S$ be a similarity in  $\rr^3$ which is  a composition of a homothety with a fixed point $c$ and a rotation around line $l\ni c$ and  let $P$ be a plane containing  line $l$. Let $B(a,r)$ be a closed ball disjoint  from $P$.  Then
$$\max\limits_{x,y\in B}\dfrac{||S(x)-x||}{||S(y)-y||}\le\dfrac{d(a,P)+r}{d(a,P)-r}\quad\quad \blacksquare $$
 
\end{lem}

\begin{lem}\label{threevecs}
Suppose $0<\al <\pi/2$,  and $\Be_1,\Be_2,\Be_3$  are such unit vectors  that the  angles  between each   $\Be_i,\Be_j$ lie  between $\pi/2-\al$
and  $\pi/2+\al$. Let $\la_1<r_i<\la_2$. Then, for  any  $\bx=(x_1, x_2, x_3)\in\rr$,
$${\la_1}(\sqrt{1-2\sin\al})\|\bx\|\le\left\|\sum\limits_{i=1}^3 x_ir_i\Be_i \right\|\le \la_2 (\sqrt{1+2\sin\al})\|\bx\|.\quad\quad\quad \blacksquare$$

\end{lem}

Using   Lemma \ref{partition} and  Lemma \ref{allaboutV}, we  find possible $\xi$-coordinates of  the point $x^\xi_k$:

\begin{lem}\label{angles}
If\ \    $S'_k(A)\cap S''_k(B)\neq\0$, then the point $x^\xi_k$ has $\xi$-coordinates
$(r_k, \phi,\be_0)$, where  $r_k=\sqrt{5} \Lip(S'_k)$,   $-\al_c< \phi<\al_c$, and
$$\al_c=\arccos\dfrac{\tan(\be-\mu/2)}{\tan(\be+\mu)}+\sqrt{5}\mu=0.295\quad\quad\blacksquare$$
\end{lem}

Using  last three  Lemmas we  show  that the map, which  assigns  to  each $\xi\in D$ the  point $x^\xi_{k+1}$, is  bi-Lipschitz:

\begin{lem}\label{Sbilip}
 For  any  k, the map  $g_k(\xi)=S^\xi_{m+1}S^{i_k}_1(z_{m-4})$ is  bi-Lipschitz with  respect  to $\xi$ on $D$.
\end{lem}

\dok.  Since $\xi,\eta\in D$, $|\Da \rho|<0.02 \rho$, $|\Da \fy|<2\mu$, $|\Da \te|<\mu/2$,  the  distances from $x_k$ to the  closest fixed point  of the  maps $H_\rho, R_\fy$ and $R_\te$ respectively are $r_k$, $r_k\sin \be_0=\dfrac{r_k}{\sqrt{5}}$ and $r_k\sqrt{\cos^2\be_0+ \sin^2\be_0\cos^2\fy}$.
Each of these  values  lies  between $r_k/\sqrt{5}$ and  $r_k$.

The displacement vectors corresponding  to each  of  these maps, have  norms 
$\dfrac{r_k\Da \rho}{ \rho}$, $\dfrac{2r_k}{\sqrt{5}}\sin\dfrac{\Da \fy}{2}$ and $2r_k\sqrt{\dfrac{5-4\sin^2\fy}{5} }\cdot \sin^2\dfrac{\Da\te}{2}$, whose
values, divided by $r_k\Da\rho$, $r_k\Da\fy$ and $r_k\Da\te$ respectively   lie  in  the interval$\left(\dfrac{0.98}{\sqrt{5}},1\right)$.

The angles  between each  two of these vectors belong to $\left(\dfrac{\pi}{2}-\al_c, \dfrac{\pi}{2}+\al_c\right)$. 
Applying  Lemma \ref{threevecs} and taking  into  account,  that $\sin\al_c<0.295$, we  have 
$${\dfrac{0.98r_k}{\sqrt{5}}}\sqrt{0.41}\|\Da\xi\|<\left\|x^\eta_k-x^\xi_k \right\|< r_k \sqrt{1.59}\|\Da\xi\|\quad\blacksquare$$

From  the  last  inequality  in the  proof  
we  obtain
$$          1.26\dfrac{\da^*_k}{r_k}    <\sqrt{\Da \rho^2+\Da \fy^2+\Da \te^2}<3.56\dfrac{\da^*_k}{r_k}  $$

The maximal displacement in the  subset $V_{m+1}$ is  reached at  the  point $z_{m+1}$ and  is less or equal  to
$ \rho\sqrt{\Da \rho^2+\Da \te^2}<3.56\dfrac{\da^*_k}{r_k} \rho<3.64\dfrac{\da^*_k}{r_k}$.

Since $\Da \te<3.56\dfrac{\da^*_k}{r_k}$,  by Lemma \ref{others}    for  any  $x\in V$,
$\|S_{m+4}^\xi(x)-S_{m+4}^\eta(x)\|<1.64\dfrac{\da^*_k}{r_k}$, so we have

\begin{lem}\label{delta}

For  any $x\in V$ and  any $i=1,...,2m$, 
$$ \|S_i^\eta(x)- S_i^\xi(x)\|<3.64\dfrac{\da^*_k}{r_k}\quad\quad\quad\quad\blacksquare$$

\end{lem}

Denote this upper bound  for the  displacement, $3.64\dfrac{\da^*_k}{r_k}$ by $\da_k$.

\bigskip
\bigskip

\subsection{ Estimates for $ F_{\xi\eta}$}

Now  we take $(i_k,j_k)\in \Sigma$, $\xi, \eta\in
D$ and estimate  the  distances between  the points  of  $S^\xi_{m+1}
S_1^{i_k}(\ga^\xi_A)$ and $S^\eta_{m+1}
S_1^{i_k}(\ga^\eta_A)$  having  the  same  addresses. 

To apply  the Proposition \ref{collage}  we    first  prove  the  following

\begin{lem}\label{bilipm+1} 
For  any  $x\in S_1^{i_k}(A)$,
$\da_k^*/1.19<\| S^\xi_{m+1}(x)-
S^\eta_{m+1}(x)\|<1.19\da_k^*$
\end{lem}
\dok

The map $F_{\xi\eta}$ is  a  composition of  a homothety  with ratio $\rho_2/\rho_1$    and  a rotation $ R_{\xi\eta}$. The map $ F_{\xi\eta}$ sends  the  point with spherical coordinates $(\rho_1,0,0)$ to the point $(\rho_2, 0, \Da\te)$, so $ R_{\xi\eta}$ sends   $(1,0,0)$ to  $(1, 0, \Da\te)$. Therefore the axis $l$ of the  rotation $ R_{\xi\eta}$ lies in the middle  bisector  plane  for  these points and the unit normal  vector to this  plane has the coordinates $\left(1,0,\dfrac{\pi+\Da\te}{2}\right)$.

\medskip

Therefore  the  set $S_{m+1}^\xi S_1^{i_k}(A)$ can be  covered  by  a  ball $W_k$
whose  radius  is $0.036 R_k$  ( where $R_k=q_{m+1}^\xi q_1^{i_k} R$) and  whose  center  has  a coordinate $(1.03 R_k,\phi, \be_0-\mu/2)$.  It  follows  from Lemma \ref{angles} that $-\al_c\le\phi\le\al_c$.

\medskip

The  maximal  angle  between the coordinate polar  axis  $l_\xi$ and  the  plane $P_l$, containing the   axis $l$  of $F_{\xi\eta}$  is $\mu/4$, therefore minimal possible  distance  between the  center  of $W_k$ and the  plane $P_l$ is greater than $0.43 R_k$.

It follows  from  the  Lemma \ref{inaball} that $\max\limits_{x,y\in W_k} \dfrac{|| F_{\xi\eta}(x)-x||}{||F_{\xi\eta}(y)-y||}\le 1.19$.

Taking $x=x_k$ and  $|| F_{\xi\eta}(x_k)-x_k||=\da_k^*$,  for  any  other  point $y\in S_{m+1}^\xi S_1^{i_k}(A)$, we have
$$\da_k^*/1.19<|| F_{\xi\eta}(y)-y||<1.19\da_k^* \quad\quad\quad\blacksquare $$

\subsection {Proof of Theorems \ref{genpos} and Proposition\ref{collage}; proof of  Theorem\ref{main}.}\label{finish}   

{\bf Proof  of  the  Theorem   \ref{genpos}.}

Since functions $\fy(x,t)$ and $\psi(x,t)$ are $\al$-H\"older
with respect to $t$, then $f(x,s,t)$  is also $\al$-H\"older with respect to $(s,t)$.

Since  $f(x,s,t)$ is  biLipschitz  with respect  to $x$,
it is clear that if $f(x_{1},s,t)=f(x_{2},s,t)$, then $x_{1}=x_{2}$,
so the equation $f(x,s,t)=0$ specifies an implicit function $x=g(s,t)$ defined  on some closed subset  $P\IN [0,1]\times [0,1]$.

Take $x_1=g(s_1,t_1)$ and $x_2=g(s_2,t_2)$ in $g(P)$. Since $f$ is  bi-Lipschitz w.r.t. $x$,
$\|f(x_2,s_1,t_1)\|\ge\|x_1-x_2\|/L$. Since $f$ is $\al$-H\"older w.r.t. $(s,t)$,
$\|f(x_2,s_1,t_1)\|\le M\|(s_1-s_2, t_1-t_2)\|^\al$.\\
Therefore $\left|x_{1}-x_{2}\right|\leq LM\left\|(s_1-s_2,t_1-t_2)\right\|^{\alpha}$.

Thus, the function $g(s,t)$ is $\al$-H\"older on the set $P$.

From $\dim_H\le 2$ we  obtain  $\dim_H(g(P))\le\dfrac{2}{\al}$.
$\blacksquare$

\bigskip

{\bf Proof  of  the  Proposition  \ref{collage}.}
For $\sa=i_1i_2i_{3}.....$ denote $\sa_k=i_{k+1}i_{k+2}....$. For $\bj=j_1...j_k$, $\hat\sa_\bj$ is an  operator sending $i_1i_2i_{3}.....$ to  $j_1...j_ki_1i_2i_{3}.....$.

\probel

By Barnsley  Collage Theorem,
$|\psi(\sa)-\fy(\sa)|\le\dfrac{\da}{1-q}$.

Let $\bi=i_1...i_k$,  $\bj=j_1...j_l$. \\Write $\psi(\sa)-\fy(\sa)=(T_{\bi}\psi(\sa_k)-S_{\bi}\psi(\sa_k))+(S_{\bi}\psi(\sa_k)-S_{\bi}\fy(\sa_k))$. 
We have $ \da_1\le\|T_{\bi}\psi(\sa_k)-S_{\bi}\psi(\sa_k)\|\le\da_2$ 
 and
$\|S_{\bi}\psi(\sa_k)-S_{\bi}\fy(\sa_k)\|\le\dfrac{q_\bi\da}{1-q}.$, which implies {\bf B1}\  \ \

The same way $\psi(\sa)-\fy(\sa)-\psi(\tau)+\fy(\tau)=$\\
$\left[(T_{\bi}\psi(\sa_k)-S_{\bi}\psi(\sa_k))-(T_{\bj}\psi(\tau_k)-S_{\bj}\psi(\tau_k))\right]+$\\
$\left[(S_{\bi}\psi(\sa_k)-S_{\bi}\fy(\sa_k))-(S_{\bj}\psi(\tau_k)-S_{\bj}\fy(\tau_k))\right]$.

The  norm  of  first  brackets lies  between $\da_1$ and $\da_2$, the  norm  of  second brackets  is no  greater  than $  \dfrac{(q_\bi+q_\bj)\da}{1-q}$,  which  gives  us {\bf B2}.
\vse

\probel

Let $\fy^\xi:[0,1]\to \ga^\xi$  be  the  linear  parametrization  of  the  zipper $\eS_\xi$ defined  in  Lemma \ref{oneforall}.  Let $I_{Ak}= T_{m+1} T_1^{i_k}(I_A)$  and $I_{Bk}=T_{m} T_{2m}^{j_k}(I_B)$ be  the  subintervals  of $I=[0,1]$,  for  which $\fy^\xi(I_{Ak})= S_{m+1} S_1^{i_k}(\ga_A)$  and $\fy^\xi(I_{Bk})= S_{m} S_{2m}^{j_k}(\ga_B)$.  Denote $\fy^\xi|_{I_{Ak}}$ by $\fy(\xi,t)$ and $\fy^\xi|_{I_{Bk}}$ by $\psi(\xi,t)$. Take $s\in I_{Ak}$, $t\in I_{Bk}$. 

\medskip

\begin{lem} \label{fxist}
 The function $f(\xi,s,t)=\fy(\xi,s)-\psi(\xi,t)$ is bi-Lipschitz with  respect  to $\xi$ for  any $s\in I_{Ak}$,   $t\in I_{Bk}$.

\end{lem}

\dok

Observe  that  if $ S_{m+1} S_1^{i_k}(A)\cap S_{m}  S_{2m}   ^{j_k}(B)\neq\0$, then
$  \dfrac{q_{m}  q_{2m}   ^{j_k}}{q_{m+1} q_1^{i_k}}<1.06$.

We  apply  the   Proposition  \ref{collage} to  $ S_{m+1} S_1^{i_k}(A)$ and  $ S_{m}  S_{2m}   ^{j_k}(B)$.

Notice  that $q_{m+4}$ is greater  than $q_{m+3} $  and $q_{m+5}$ and  the  same  is  true for their  symmetric  counterparts,  so  we  use  $q_{m+4}$ (which  is  equal to $q_{m-4}$) for both $S_{m+1} S_1^{i_k}(A)$  and     $S_{m}  S_{2m}   ^{j_k}(B)$.

By Proposition  \ref{collage}, for  any    $s\in I_{Ak}$,   $t\in I_{Bk}$,
$$\dfrac{\da_k^*}{1.19}-\dfrac{q_{m+4}\da_k}{1-q}(q_{m+1} q_1^{i_k}+q_{m}  q_{2m}   ^{j_k})\le |\fy(\xi,s)-\fy(\eta,s)-\psi(\xi,t)+\psi(\eta,t)|\le$$  $$1.19\da_k^* +\dfrac{q_{m+4}\da_k}{1-q}(q_{m+1} q_1^{i_k}+q_{m}  q_{2m}   ^{j_k})$$
Evaluating
$\dfrac{q_{m+4}}{1-q}(q_{m+1} q_1^{i_k}+q_{m}  q_{2m}   ^{j_k})\da_k<0.02\cdot2.06q_{m+1} \cdot3.64\da^*<0.03\da_k^*$, we  get
$$0.8\da_k^*< |\fy(\xi,s)-\fy(\eta,s)-\psi(\xi,t)+\psi(\eta,t)|< 1.22\da_k^* $$
Since for  the  point $x_k$, the  function          $\Fy(\xi,\eta)= F_{\xi\eta}(x_k)-x_k$ is  bi-Lipschitz  with respect  to  $\xi$ and $\eta$, and  since $\| F_{\xi\eta}(x_k)-x\|=\da_k^*$, the  last  inequality  shows  that  the  function $f(\xi,s,t)=\fy(\xi,s)-\psi(\xi,t)$ is bi-Lipschitz with  respect  to $\xi$. \vse

\section{Addendum: Dense   groups  of   second type and  their  generators.}

We will remind  several  facts  about dense 2- generator  subgroups in $\bbc^*$. See [3] for  details.

As it follows  from Kronecker's Theorem

\begin{prop}\label{mulgr} Let $u,v \in \mathbb{C}$,  $Im \dfrac{u}{v} \neq 0$, $\alpha u + \beta v=1$, $\alpha,\beta \in \mathbb{R}$, $\xi=e^{2\pi i u}$,    $\eta=e^{2\pi i v}$.\\
A group $G=\langle \xi,\eta,\cdot\rangle $ is  dense in $\mathbb{C}$
 iff
 for any integers $k,l,m$,
$$ k\alpha+l\beta+m =0  \mbox{\rm \ \  implies \ \ }  k=l=m=0. \ \ \ \  \blacksquare $$
\end{prop}

Then $G=\langle \xi,\eta,\cdot\rangle $ is called a  {\em dense 2-generator multiplicative  group}  in $\mathbb{C}$.

\bigskip

Notice that if  $G=\langle \xi,\eta,\cdot\rangle $ is dense in $\mathbb{C}$, the  formula $\psi(\xi^m\eta^n)=\left(\dfrac{\xi}{|\xi|}\right)^m$  defines
  a homomorphism $\psi$ of  a group $G=\langle \xi,\eta,\cdot\rangle $ to the unit circle
$S^1\subset \mathbb{C}$.
Put $$H_G= \bigcap\limits_{\varepsilon>0}\overline{\psi(B(1,\varepsilon)\cap G)}$$

In other  words, $H_G$ is the  set of limit points of all those  sequences $\{e^{i n_k\arg(\xi)}\}$, for  which   $\{n_k\}$ are the first  coordinates    of such   sequence $\{(n_k,m_k)\}$, that $ \lim\limits_{k\to\infty} \xi^{n_k}\eta^{m_k}=1$.

The set $H_G$ is a closed topological subgroup of the unit  circle $S^1$, so it is either finite cyclic or  infinite.

\begin{dfn}
A dense 2-generator subgroup $G$ is called  the group of   first  type, if $H_G$ is finite, and  the group of  second type, if $H_G=S^1$.
\end{dfn}

 So $G$ is of   second type iff for  some $\alpha\notin \mathbb{Q}$,  $e^{2\pi i \alpha}\in H_G$.

Therefore, if $\xi,\eta$ are the generators of a group of  second type, then for any rational $p,q$, the numbers $\xi^p,\eta^q$ also generate a group of the second type. This implies

\begin{prop}\label{dense in C2}
The  set of  pairs  $\xi,\eta$ of  generators  of the groups of   second type is  dense  in $\bbc^2$. \vse\end{prop}

The  groups of   second type  have  a significant  geometric  property:

\begin{thm}
If  the  group $G=\langle \xi,\eta,\cdot\rangle$, $\xi=r e^{i\alpha},\eta=R e^{i\beta}$ is  of  second type, then for  any   $z_1,z_2\in \mathbb{C}\setminus\{0\}$ there is such sequence $\{(n_k,m_k)\}$  that $\lim\limits_{k\to \infty}   \dfrac{z_1\xi^{n_k}}{z_2\eta^{m_k}}=1$,  $\lim\limits_{k\to\infty} e^{i n_k\alpha}=e^{-i \arg(z_1)}$,  $\lim\limits_{k\to\infty} e^{i n_k\beta}=e^{-i \arg(z_2)}$.
\end{thm}

\end{document}